# TENSOR PRODUCTS AND RELATION QUANTALES

MARCEL ERNÉ AND JORGE PICADO

ABSTRACT. A classical tensor product $A \otimes B$ of complete lattices $A$ and $B$, consisting of all down-sets in $A \times B$ that are join-closed in either coordinate, is isomorphic to the complete lattice $Gal(A, B)$ of Galois maps from $A$ to $B$, turning arbitrary joins into meets. We introduce more general kinds of tensor products for closure spaces and for posets. They have the expected universal property for bimorphisms (separately continuous maps or maps preserving restricted joins in the two components) into complete lattices. The appropriate ingredient for quantale constructions is here distributivity at the bottom, a generalization of pseudocomplementedness. We show that the truncated tensor product of a complete lattice $B$ with itself becomes a quantale with the closure of the relation product as multiplication iff $B$ is pseudocomplemented, and the tensor product has a unit element iff $B$ is atomistic. The pseudocomplemented complete lattices form a semicategory in which the hom-set between two objects is their tensor product. The largest subcategory of that semicategory has as objects the atomic boolean complete lattices, which is equivalent to the category of sets and relations. More general results are obtained for closure spaces and posets.

## 1. INTRODUCTION

Tensor products have their place in algebra, (point-free) topology, order theory, category theory and other mathematical disciplines. In the realm of ordered sets, they are intimately related to the concept of Galois connections (see e.g. [5, 11, 19, 30, 34]; for historical background, refer to [14, 16, 18]). In the present paper, we show how such tensor products give rise to certain quantales whose members are specific relations between complete lattices, partially ordered sets (posets) or closure spaces.

Before focussing on tensor products, let us recall briefly the fundamental notions in the theory of Galois connections. Given two posets $A$ and $B$, let $Ant(A, B)$ denote the pointwise ordered set of all *antitone*, i.e. order-reversing maps from $A$ to $B$, and $Gal(A, B)$ the subposet of all *Galois maps*, i.e. maps from $A$ to $B$ such that the preimage of any principal filter is a principal ideal [14, 28, 34]. If $A$ and $B$ are complete lattices, $Ant(A, B)$ is a

*Date*: December 16, 2016.
2010 *Mathematics Subject Classification.* Primary: 06F05. Secondary: 08A30, 03G05, 16D30, 20M12, 20N02, 16N80, 17A65.
The second author acknowledges support from CMUC (UID/MAT/00324/2013 funded by the Portuguese Government through FCT/MCTES and co-funded by the European RDF through Partnership Agreement PT2020) and grants MTM2015-63608-P (Ministry of Economy and Competitiveness of Spain) and IT974-16 (Basque Government).

complete lattice, too, and $Gal(A,B)$ is the complete lattice of all $f\colon A\to B$ satisfying $f(\bigvee X)=\bigwedge f[X]$ for each $X\subseteq A$. Galois maps are closely tied to *Galois connections*; these are dual adjunctions between posets $A$ and $B$, that is, pairs $(f,g)$ of maps $f\colon A\to B$ and $g\colon B\to A$ such that

$$x\leq g(y)\Leftrightarrow y\leq f(x)\ \text{ for all } x\in A \text{ and } y\in B,$$

or equivalently, pairs of maps $f\in Ant(A,B)$ and $g\in Ant(B,A)$ with

$$x\leq g(f(x)) \text{ for all } x\in A\ \text{ and }\ y\leq f(g(y)) \text{ for all } y\in B.$$

Either partner in a Galois connection determines the other by the formula

$$g(y)=f^*(y)=\max\{x\in A : f(x)\geq y\},$$

and the Galois maps are nothing but the partners of Galois connections. Clearly, $(f,g)$ is a Galois connection iff $(g,f)$ is one, and consequently, $Gal(A,B)\simeq Gal(B,A)$. Both composites of the partners of a Galois connection are closure operations, and their ranges are dually isomorphic.

Three tensor products of posets $A$ and $B$ are defined as follows: $A\otimes_r B$ denotes the collection of all *right tensors*, i.e. down-sets $T$ in $A\times B$ such that $xT=\{y\in B : (x,y)\in T\}$ is a principal ideal of $B$ for each $x\in A$. The system $A_\ell\otimes B$ of *left tensors* is defined in the opposite manner, and the *tensor product* $A_\ell\otimes_r B$ consists of all (*two-sided*, i.e. left and right) *tensors*; if $A$ and $B$ are complete, it is denoted by $A\otimes B$; in that case, a down-set $T$ in $A\times B$ is a right tensor iff $\{x\}\times Y\subseteq T$ implies $(x,\bigvee Y)\in T$, a left tensor iff $X\times\{y\}\subseteq T$ implies $(\bigvee X,y)\in T$, and a tensor (or *G-ideal* [31]) iff $X\times Y\subseteq T$ implies $(\bigvee X,\bigvee Y)\in T$ (see Shmuely [34] for alternative characterizations).

A bijective connection between posets of antitone maps and tensor products of posets $A,B$ is provided by the assignments

$$f\mapsto T_f=\{(x,y)\in A\times B : f(x)\geq y\}\ \text{ and}$$
$$T\mapsto f_T\text{ with } f_T\colon A\to B, x\mapsto \max xT.$$

Indeed, these maps are mutually inverse isomorphisms between $A\otimes_r B$ and $Ant(A,B)$, and they induce isomorphisms between $A_\ell\otimes_r B$ and $Gal(A,B)$.

If some poset $B$ has a least element $0=0_B$, we may build the "truncated" poset $\check{B}=B\smallsetminus\{0_B\}$. Now, given complete lattices $A,B,C$ and $f\in Ant(A,B)$, $g\in Ant(B,C)$, define $g\odot f\colon A\to C$ by

$$g\odot f(a)=\bigvee\{c\in C : (a,c)\in E_{f,g}\},$$

where $E_{f,g}$ denotes the tensor generated by the set

$$T_{f,g}=\{(x,z)\in A\times C : \exists\, y\in\check{B} : f(x)\geq y \text{ and } g(y)\geq z\}$$
$$(\,=\{(x,z)\in A\times C : f(x)\wedge g^*(z)>0_B\} \text{ if } g \text{ is a Galois map}).$$

Picado showed in Proposition 3.1 of [31] that the so-defined $g\odot f$ is in fact a Galois map from $A$ to $C$. This gives a way of composing antitone maps and Galois maps or connections, so that the composed map is again antitone, which almost never would happen with the usual composition of maps. In certain cases, the alternate composition $\odot$ appears somewhat mysterious.

**Example 1.1.** If $\mathbb{I}$ denotes the real unit interval $[0,1]$ with the usual order, the composite $g \odot f$ of $f, g \in Ant(\mathbb{I}, \mathbb{I})$ is always a step function! Explicitly,

$$T_{f,g} = \{(x,z) : \exists\, y > 0 \; (f(x) \geq y, \; g(y) \geq z)\}$$
$$= \{(x,z) : f(x) > 0, \; \exists\, y > 0 \; (g(y) \geq z)\},$$
$$E_{f,g} = \{(a,c) : a \leq r = \bigvee\{x : f(x) > 0\}, \; c \leq s = \bigvee\{g(y) : y > 0\}\} \cup \overline{\emptyset},$$

where $\overline{\emptyset} = (\{0\} \times \mathbb{I}) \cup (\mathbb{I} \times \{0\})$. Therefore,

$g \odot f(0) = 1$, $g \odot f(a) = s$ if $0 < a \leq r$, and $g \odot f(a) = 0$ otherwise.

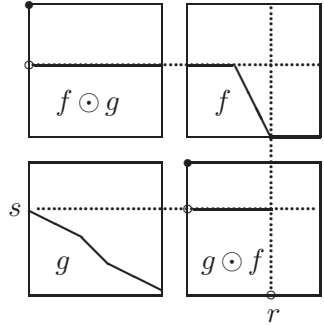

In particular, the new composition $\odot$ of any two involutions (that is, antitone bijections) of $\mathbb{I}$ yields the constant function 1. Apparently, the previous arguments work in any complete chain instead of $\mathbb{I}$, and much more generally, in any complete lattice with a meet-irreducible least element 0.

In [31] it is proved that $\odot$ makes $Gal(B, B)$ a quantale whenever $B$ is a frame (locale). Hence, with any frame $B$, there is associated not only a quantale of (isotone, i.e. order preserving) residuated maps [6, 18], but also a quantale of (antitone!) Galois maps. One of our main goals is to characterize those complete lattices $B$ for which $Gal(B, B)$ or $Ant(B, B)$, respectively, together with the multiplication $\odot$ becomes a quantale. In most cases, it will be technically more comfortable to work with tensor products than with $Gal(A, B)$ or $Ant(A, B)$.

After a brief summary of some central notions and facts in quantale theory we introduce binary tensor products for closure spaces and establish their fundamental properties. This allows us to extend the theory of tensor products for complete lattices in diverse directions, in particular, from complete lattices to arbitrary posets. The usual trick is here to replace joins with cuts (see [10, 11, 16] and Section 4), and then, in a more courageous step, cuts by arbitrary closed sets in closure spaces. A *tensor* between closure spaces $A$ and $B$ is a subset $T$ of $A \times B$ such that all "slices" $xT$ and $Ty$ are closed, and the tensor product $A \otimes B$ is the closure system of all such tensors.

Any *augmented poset* $A_{\mathcal{X}} = (A, \mathcal{X})$ (where $A$ is a poset and $\mathcal{X}$ a collection of subsets of $A$) may be interpreted as a closure space, by considering the closure system of all $\mathcal{X}$-ideals, i.e. down-sets $I$ containing the cut closure $\Delta X = \bigcap\{\downarrow y : X \subseteq \downarrow y\}$ whenever $X \in \mathcal{X}$ and $X \subseteq I$. Then, the tensor product $A_{\mathcal{X}} \otimes B_{\mathcal{Y}}$ of two augmented posets consist of all $\mathcal{X}\mathcal{Y}$-*ideals* or $\mathcal{X}\mathcal{Y}$-*tensors*, i.e. down-sets $T$ in $A \times B$ such that for all $X \in \mathcal{X}$ and $Y \in \mathcal{Y}$, $X \times Y \subseteq T$ implies $\Delta X \times \Delta Y \subseteq T$ (resp. $(\bigvee X, \bigvee Y) \in T$ if the involved joins exist). If $A$ and $B$ are complete lattices then the tensor product $A_{\mathcal{X}} \otimes B = A_{\mathcal{X}} \otimes B_{\mathcal{P}B}$ is isomorphic to the complete lattice $Ant_{\mathcal{X}}(A, B)$ of maps $f \colon A \to B$ satisfying $f(\bigvee X) = \bigwedge f[X]$ for all $X \in \mathcal{X}$. Furthermore, our tensor products have the expected universal property with respect to the appropriate bimorphisms.

A considerable simplification is achieved by passing to *truncated tensor products* $A_{\mathcal{X}} \check{\otimes} B_{\mathcal{Y}}$, cutting off the least tensor from all tensors. Their elements are down-sets in the direct product $\check{A} \times \check{B}$ (with $\check{A} = A \smallsetminus \overline{\emptyset}$ and $\check{B} = B \smallsetminus \overline{\emptyset}$) such that the conditions in the two coordinates hold for *nonempty* "rectangles". One advantage of that reduction is that the pure tensors $a \otimes b$ have no longer the rather complicated form $\overline{(a,b)} \cup (\overline{\emptyset} \times B) \cup (A \times \overline{\emptyset})$ but become simply point closures, resp. principal ideals. Another, and more important, advantage is that now the quantale constructions are much easier, since the tensor multiplication corresponding to $\odot$ is obtained by forming the (right) tensor closure of the usual relation product, and then the order isomorphism between $A_{\ell} \otimes_r B$ and $Gal(A, B)$, resp. between $A \otimes_r B$ and $Ant(A, B)$, also transports the multiplication. A main result will be that for any complete lattice $B$, the truncated tensor product $B \check{\otimes} B$, resp. the isomorphic tensor product $B \otimes B \simeq Gal(B, B)$, becomes a quantale iff $B$ is pseudocomplemented, and a unital quantale iff $B$ is an atomic boolean complete lattice.

Our constructions also provide a *semicategory* (missing identity morphisms) of pseudocomplemented complete lattices together with the (truncated) tensors or antitone maps, respectively, as morphisms. In that semicategory, the atomic boolean complete lattices (isomorphic copies of power set lattices) form the greatest subcategory, and the latter is equivalent to the category of sets and relations as morphisms. Similar results are obtained for augmented posets and for closure spaces instead of complete lattices. Crucial is here the observation that a closure system is pseudocomplemented iff all *polars* $x^{\perp} = \{y : \overline{x} \cap \overline{y} = \overline{\emptyset}\}$ are closed.

In the case of a frame / locale, our quantale constructions via Galois connections or tensor products have important applications to the point-free treatment of uniform structures (see Ferreira and Picado [20]). For a comprehensive discussion of frames resp. locales see Picado and Pultr [32], and for more background concerning quantales Rosenthal [33].

For categorical topics refer to Adámek, Herrlich and Strecker [1]. The diverse relation products discussed in this paper fit, of course, into the general category-theoretical framework of relations and their multiplication (see [1], Ferreira and Picado [20, Section 2], Freyd and Scedrov [21]).

The main results in this paper have been presented by the first author at the 71st Workshop on General Algebra (AAA), Bedlewo 2006 [15] but were never published in a journal until now.

## 2. Prequantales, quantales, prenuclei and nuclei

A *prequantale* is a complete lattice $Q$ with a multiplication that distributes over arbitrary joins from both sides: for all $x$ in $Q$ and subsets $Y$ of $Q$,

$$x \cdot \bigvee Y = \bigvee\{x \cdot y : y \in Y\} \text{ and } \bigvee Y \cdot x = \bigvee\{y \cdot x : y \in Y\}.$$

A *quantale* is a prequantale with associative multiplication, and a *locale* or *frame* is a quantale whose multiplication is the binary meet. An equivalent definition of quantales characterizes them as *complete residuated semigroups* (Birkhoff [5], Blyth and Janowitz [6]) with residuation operations satisfying

$$y \leq z \mathbin{.\!\cdot} x = x \rightarrow z \Leftrightarrow x \cdot y \leq z \Leftrightarrow x \leq z \mathbin{\cdot\!.} y = z \leftarrow y.$$

A *subquantale* of a quantale $Q$ is closed under joins and the multiplication, and a *quantale homomorphism* is a map between quantales that preserves arbitrary joins and the multiplication. The example below is taken from [18].

**Example 2.1.** Let $(X, \cdot)$ be a *partial semigroup*, that is, a set equipped with a partial binary operation such that whenever $r \cdot s$ and $s \cdot t$ are defined, then so are $(r \cdot s) \cdot t$ and $r \cdot (s \cdot t)$, and these products are equal. Any such partial semigroup gives rise to a quantale whose ground set is the power set $\mathcal{P} X$ and whose multiplication and residuals are given by

$$R \odot S = \{r \cdot s : r \in R, \, s \in S, \text{ and } r \cdot s \text{ is defined}\},$$
$$R \rightarrow T = \{x \in X : \text{ for all } r \in R, \text{ if } r \cdot x \text{ is defined then } r \cdot x \in T\},$$
$$T \leftarrow S = \{x \in X : \text{ for all } s \in S, \text{ if } x \cdot s \text{ is defined then } x \cdot s \in T\}.$$

Specifically, taking the partial semigroup operation on a cartesian square $B \times B$ given by $(a,b) \cdot (d,c) = (a,c)$ if $b = d$, the power set $\mathcal{P}(B \times B)$ becomes a quantale with unions as joins and the usual relation product

$$S \circ R = R \cdot S = \{(a,c) : (a,b) \in R \text{ and } (b,c) \in S \text{ for some } b\}.$$

A *preclosure operation* on a poset $B$ is an extensive (increasing) and isotone (order preserving) map from $B$ to $B$, and a *closure operation* is an idempotent preclosure operation; the fixpoint set

$$B_j = \{x \in B : j(x) = x\}$$

of a preclosure operation $j$ is a *closure domain*, and the corresponding closure operation $\overline{j}$ is the (pointwise) least closure operation above $j$. The poset of closure operations is dually isomorphic to that of closure domains (see e.g. [16]). As in [31], we mean by a *prenucleus* on a (pre)quantale $Q$ a preclosure operation $j \colon Q \to Q$ such that

$$x \cdot j(y) \vee j(x) \cdot y \leq j(x \cdot y) \text{ for all } x, y \in Q.$$

A (*quantic* or *quantalic*) *nucleus* on $Q$ is an idempotent prenucleus, or equivalently, a closure operation $j$ on $Q$ satisfying

$$j(x) \cdot j(y) \leq j(x \cdot y) \text{ for all } x, y \in Q.$$

It is well known and easy to see that the surjective corestrictions of nuclei on $Q$ form a set of representatives for the surjective quantale homomorphisms $f \colon Q \to Q'$: such an $f$ has an upper (right) adjoint $g$; then $k = g \circ f$ is a nucleus on $Q$, and $i \colon Q_k \to Q'$, $x \mapsto f(x)$ is an isomorphism with $f = i \circ k$.

The fixpoint sets of (pre)nuclei are characterized in Proposition 2.1 below (see Niefield and Rosenthal [29, 33] and Picado [31]). Extensions to prenuclei and nuclei on prequantales or *cm-lattices*, where the multiplication is isotone but not necessarily distributive over joins, are possible (Erné [12, 13]).

**Proposition 2.1.** *For a quantale $Q$ and $S \subseteq Q$, the following are equivalent:*

(a) $S$ is the fixpoint set $Q_j$ of some prenucleus $j$ on $Q$.
(b) $S$ is the range $Q_k$ of some nucleus $k$ on $Q$.
(c) $S$ is closed under arbitrary meets and residuation in $Q$, that is, $q \in Q$ and $s \in S$ imply $q \to s \in S$ and $s \leftarrow q \in S$.
(d) $S$ is closed under arbitrary meets in $Q$ and is a quantale with respect to the multiplication $x \cdot_S y = \bigwedge \{s \in S : x \cdot y \leq s\}$.

*If these conditions hold then $k = \overline{j}$ and $x \cdot_S y = k(x \cdot y)$. There is a bijective correspondence between nuclei and subsets $S$ satisfying these conditions.*

We call a subset $S$ enjoying the above properties a *quantic quotient* of $Q$. (Warning: a quantic quotient of a locale is usually referred to as a *sublocale*, because the locale morphisms go in the opposite direction, whereas a subquantale of a locale resp. frame is a *subframe*!)

In later sections, we shall need a related result (cf. [31, Lemma 3.2]):

**Lemma 2.1.** *Let $i\colon I \to I$, $j\colon J \to J$, $k\colon K \to K$ be preclosure operations on complete lattices, and suppose a map $m\colon I \times J \to K, (x,y) \mapsto x \cdot y$ satisfies*

$$\bigvee X \cdot \bigvee Y = \bigvee (X \cdot Y) \ (= \bigvee \{x \cdot y : x \in X, y \in Y\}) \text{ for all } X \subseteq I, Y \subseteq J.$$

*If*   $i(x) \cdot y \vee x \cdot j(y) \leq \overline{k}(x \cdot y)$   *for all $x \in I$, $y \in J$*
*then*   $\overline{k}(\overline{i}(x) \cdot \overline{j}(y)) = \overline{k}(x \cdot y)$   *for all $x \in I$, $y \in J$.*

*Proof.* Let $E$ denote the pointwise ordered set of all preclosure operations on $I$, and consider the subset

$$F = \{f \in E : \forall\, x \in I \ \forall\, y \in J \ (f(x) \cdot y \leq \overline{k}(x \cdot y))\}.$$

By hypothesis, we have $i \in F$ and, as $F$ is closed under pointwise suprema (by the distribution law), $i \leq s = \bigvee F \in F$. Thus, $s$ is the greatest element of $F$. But $s \leq s^2 \in E$ and $s^2(x) \cdot y \leq \overline{k}(s(x) \cdot y) \leq \overline{k}(\overline{k}(x \cdot y)) = \overline{k}(x \cdot y)$ imply $s^2 \in F$. It follows that $s = s^2$ is a closure operation, and consequently $\overline{i} \leq s$, whence $\overline{i}(x) \cdot y \leq \overline{k}(x \cdot y)$. Analogously, we get $x \cdot \overline{j}(y) \leq \overline{k}(x \cdot y)$, and finally $\overline{i}(x) \cdot \overline{j}(y) \leq \overline{k}(\overline{i}(x) \cdot y) \leq \overline{k}(\overline{k}(x \cdot y)) = \overline{k}(x \cdot y)$, hence $\overline{k}(\overline{i}(x) \cdot \overline{j}(y)) \leq \overline{k}(x \cdot y)$. The reverse inequality follows from $x \cdot y \leq \overline{i}(x) \cdot \overline{j}(y)$. □

Similar arguments show that many other (in)equalities involving $i, j, k$ may be transferred to the corresponding (in)equalities for $\overline{i}, \overline{j}, \overline{k}$.

## 3. Tensor products of closure spaces

By slight abuse of language, one calls a closure (resp. preclosure) operation on a power set $\mathcal{B}=\mathcal{P}E$ a *closure operator* (resp. *preclosure operator*) on $E$, and a collection of subsets closed under arbitrary intersections (with $\bigcap \emptyset = E$) a *closure system* on $E$. The pair $A = (E,\mathcal{C})$ is then referred to as a *closure space*, and one writes $\mathcal{U}A$ for $E$ and $\mathcal{C}A$ for $\mathcal{C}$. For each preclosure operator $c$ on $E$, the set $\mathcal{B}_c$ of its fixpoints is a closure system. In the opposite direction, one assigns to any subset $\mathcal{X}$ of $\mathcal{B}$ the closure operator $c_{\mathcal{X}} \colon \mathcal{B} \to \mathcal{B}$ with

$$c_{\mathcal{X}}(Y) = c_{\mathcal{X}} Y = \overline{Y} = \bigcap \{X \in \mathcal{X} : Y \subseteq X\}.$$

The assignments $c \mapsto \mathcal{B}_c$ and $\mathcal{X} \mapsto c_{\mathcal{X}}$ constitute a Galois connection that induces a dual isomorphism between the complete lattice of all closure operators and that of all closure systems on $E$ [5, 16, 30].

Open sets, subspaces, and continuous maps between closure spaces are defined as in the topological case. Given any preclosure operator $c$ on $E$, a map $f$ is continuous between the closure spaces $(E,\mathcal{B}_c)$ and $(F,\mathcal{C})$ iff the inclusion $f[c(Y)] \subseteq c_{\mathcal{C}} f[Y]$ holds for all $Y \subseteq E$. A closure space is $T_0$ iff distinct singletons have different closures. For any complete lattice $B$ and any subset $X$ of $B$, the pair $(X, \{\{x \in X : x \leq b\} : b \in B\})$ is a $T_0$ closure space, and *all* $T_0$ closure spaces arise in that fashion [16, 30]. Every closure space has an obvious $T_0$ reflection (obtained by identifying points with the same closure) whose closure system is isomorphic to the original one. Thereby, great parts of the theory of closure spaces may be reduced to complete lattices (cf. [9]).

The empty set often causes some troubles in the subsequent considerations. Therefore, we call a closure space $A$ *unbounded* if $\emptyset$ is closed, and *bounded* otherwise, because then any element in the least closed set $\bot_A = \overline{\emptyset}$ is a lower bound of the space in the *specialization order* given by $x \leq y \Leftrightarrow x \in \overline{y} = \overline{\{y\}}$. The closure space is *uniquely bounded* if the closure of $\emptyset$ is a singleton. As we saw above, all complete lattices give rise to uniquely bounded closure spaces whose closed sets are the principal ideals. The restriction to bounded or to unbounded closure spaces does not cause a severe loss of generality: adding a point $0$ and taking the sets $C \cup \{0\}$ instead of the original closed sets $C$, one obtains a uniquely bounded reflection; on the other hand, passing from a closure space $A$ to the subspace $\breve{A}$ on $\mathcal{U}A \smallsetminus \bot_A$ yields an unbounded coreflection; and all three spaces have isomorphic closure systems. Let us note (cf. [9, 16]):

**Proposition 3.1.** *The construct of complete lattices (regarded as closure spaces) and join-preserving maps is a full epireflective subconstruct of the construct of closure spaces: for any closure space $A = (E,\mathcal{C})$, the map*

$$\eta_A \colon A \to \mathcal{C}, \ x \mapsto \overline{x}$$

*is continuous, and for any continuous map $f$ from $A$ into a complete lattice $B$, there is a unique join-preserving map $f^{\vee} \colon \mathcal{C} \to B$ with $f = f^{\vee} \circ \eta_A$.*

The product of two closure spaces $A$ and $B$ is the cartesian product $\mathcal{U}A \times \mathcal{U}B$ of the underlying sets, having the "rectangles" $F \times G$ with $F \in \mathcal{C}A$, $G \in \mathcal{C}B$ as closed sets. Notice that products of topological spaces have much more closed sets! In the product $A \times B$ of closure spaces, one has the equation $\overline{X \times Y} = \overline{X} \times \overline{Y}$ if $X \neq \emptyset \neq Y$, which however may fail if $X$ or $Y$ are empty.

Now, we define the *tensor product* $A \otimes B$ of $A$ and $B$ to be the closure system of all *tensors*, i.e. subsets $T$ of the underlying set of $A \times B$ such that for each $x$ in $A$, the set $xT = \{y : (x,y) \in T\}$ is closed in $B$, for each $y$ in $B$, the set $Ty = \{x : (x,y) \in T\}$ is closed in $A$.

Equivalently, $X \times Y \subseteq T$ implies $\overline{X} \times \overline{Y} \subseteq T$. The closure space

$$A \overline{\times} B = (\mathcal{U}A \times \mathcal{U}B, A \otimes B)$$

is referred to as the *tensorial product* of $A$ and $B$. The least element of $A \otimes B$ is

$$\overline{\emptyset} = (\bot_A \times \mathcal{U}B) \cup (\mathcal{U}A \times \bot_B),$$

where the symbol $\overline{\emptyset}$ refers to $A \overline{\times} B$, while $\bot_A = \overline{\emptyset}$ refers to $A$ and $\bot_B = \overline{\emptyset}$ to $B$. In $A \overline{\times} B$, the following product equation is always true:

$$\overline{X \times Y} = (\overline{X} \times \overline{Y}) \cup \overline{\emptyset}.$$

Any tensor $T \in A \otimes B$ is the join (in fact, the union) of the *pure tensors* $x \otimes y$ with $(x,y) \in T$, where $x \otimes y = (\overline{x} \times \overline{y}) \cup \overline{\emptyset}$ is the least tensor containing $(x,y)$.

Given closure spaces $A, B, C$ and a map $f : A \times B \to C$, we have maps

$$_xf : B \to C, \ y \mapsto f(x,y) \ (x \in \mathcal{U}A),$$
$$f_y : A \to C, \ x \mapsto f(x,y) \ (y \in \mathcal{U}B).$$

The map $f$ is *separately continuous* if all $_xf$ and all $f_y$ are continuous.

**Lemma 3.1.** (1) *A map $f : A \times B \to C$ is separately continuous iff it is continuous as a map from $A \overline{\times} B$ to $C$.*

(2) *Nonempty closure spaces $A$ and $B$ are unbounded iff every continuous map $f : A \times B \to C$ is separately continuous.*

*Proof.* (1) follows from the identities $_xf^{-1}[Z] = xT$ and $f_y^{-1}[Z] = Ty$ for $Z \subseteq C$ and $T = f^{-1}[Z]$.

(2) Suppose $f : A \times B \to C$ is continuous. Then the preimage of any closed set of $C$ under $f$ is a closed rectangle $T = F \times G$. If $A$ and $B$ are unbounded then $T$ is closed in $A \overline{\times} B$, since each $xT$ is $G$ or $\emptyset$, and each $Ty$ is $F$ or $\emptyset$. Now, (1) applies to show that $f$ is separately continuous.

On the other hand, if, say, $A$ is bounded then $\emptyset = \emptyset \times \mathcal{U}B$ is closed in $A \times B$, but not in $A \overline{\times} B$, since $\overline{\emptyset} \times \mathcal{U}B \neq \emptyset$. Hence, the identity map on $A \times B$ is continuous but not separately continuous, by (1). □

The coincidence $A \times B = A \overline{\times} B$ is rather rare: it fails, for instance, whenever $A$ and $B$ are non-singleton closure spaces whose finite subsets are closed.

By earlier remarks, the tensor product $A \otimes B$, being a complete lattice, may be regarded as a specific closure space, with the principal ideals $\downarrow \{T\}$ as closed sets. Now, from Proposition 3.1 and Lemma 3.1 (1), we infer:

**Theorem 3.1.** *For arbitrary closure spaces $A, B$, the pure tensor insertion*

$$\otimes \colon A \times B \to A \otimes B, \ (x,y) \mapsto x \otimes y$$

*is universal among all separately continuous maps from $A \times B$ to complete lattices: a map $f$ from $A \times B$ to a complete lattice $C$ is separately continuous iff there is a unique join-preserving $f^\vee \colon A \otimes B \to C$ with $f(x,y) = f^\vee(x \otimes y)$. This map is given by $f^\vee(T) = \bigvee f[T]$.*

**Corollary 3.1.** *Tensor products of closure spaces satisfy the commutative law $A \otimes B \simeq B \otimes A$ and the associative law $(A \otimes B) \otimes C \simeq A \otimes (B \otimes C)$.*

Fortunately, tensor products of closure spaces may be reduced to those of complete lattices by the next result, which resembles a similar theorem for concept lattices (cf. [11, 36]) but refers to a different tensor product (cf. [10]):

**Theorem 3.2.** *For any two closure spaces $A$ and $B$, the complete lattices $A \otimes B$ and $\mathcal{C}A \otimes \mathcal{C}B$ are isomorphic, via the bijection*

$$h \colon A \otimes B \to \mathcal{C}A \otimes \mathcal{C}B, \ T \mapsto \{(X,Y) \in \mathcal{C}A \times \mathcal{C}B : X \times Y \subseteq T\}.$$

*Proof.* Concerning well-definedness of $h$, notice that for any $X \in \mathcal{C}A$ and any $\mathcal{Y} \subseteq X\,h(T) = \{Y \in \mathcal{C}B : (X,Y) \in h(T)\}$, one has $X \times \bigcup \mathcal{Y} \subseteq T$, hence $X \times \overline{\bigcup \mathcal{Y}} \subseteq T$, and so $\bigvee \mathcal{Y} = \overline{\bigcup \mathcal{Y}} \in X\,h(T)$; and similarly for the other side.

$h$ is an order embedding: if $S, T \in A \otimes B$ satisfy $S \subseteq T$ then $(X,Y) \in h(S)$ implies $X \times Y \subseteq S \subseteq T$, hence $(X,Y) \in h(T)$. Conversely, if $h(S) \subseteq h(T)$ and $(x,y) \in S$ then $\overline{x} \times \overline{y} \subseteq S$, $(\overline{x}, \overline{y}) \in h(S) \subseteq h(T)$, and so $(x,y) \in \overline{x} \times \overline{y} \subseteq T$.

$h$ is onto: given $\mathcal{T} \in \mathcal{C}A \otimes \mathcal{C}B$, put $T = \{(x,y) \in \mathcal{U}A \times \mathcal{U}B : h(x \otimes y) \subseteq \mathcal{T}\}$. Then we get $h(T) \subseteq \mathcal{T}$. Indeed, $(X,Y) \in h(T)$ implies $X \times Y \subseteq T$ and so $h(x \otimes y) \subseteq \mathcal{T}$ for all $(x,y) \in X \times Y$; any such $(x,y)$ satisfies $\overline{x} \times \overline{y} \subseteq x \otimes y$, hence $(\overline{x}, \overline{y}) \in h(x \otimes y) \subseteq \mathcal{T}$. Now, $\mathcal{T} \in A \otimes B$ yields $(X,Y) \in \mathcal{T}$, as $X$ is the join of the point closures $\overline{x}$ with $x \in X$, and analogously for $Y$.

Conversely, for $(F,G) \in \mathcal{T}$ and each $(x,y) \in F \times G$, we prove $h(x \otimes y) \subseteq \mathcal{T}$: indeed, $(X,Y) \in h(x \otimes y)$ implies $X \times Y \subseteq x \otimes y \subseteq (F \times G) \cup \overline{\emptyset}$ (as $(x,y)$ lies in $(F \times G) \cup \overline{\emptyset} \in A \otimes B$); hence $(X \smallsetminus \bot_A) \times (Y \smallsetminus \bot_B) = (X \times Y) \smallsetminus \overline{\emptyset} \subseteq F \times G$; if both $X \smallsetminus \bot_A$ and $Y \smallsetminus \bot_B$ are nonempty, it follows that $X \smallsetminus \bot_A \subseteq F$ and $Y \smallsetminus \bot_B \subseteq G$; consequently $(X,Y) \leq (F,G)$ in $\mathcal{C}A \otimes \mathcal{C}B$, and so $(X,Y) \in \mathcal{T}$, as $\mathcal{T}$ is a down-set (see Section 4) containing $(F,G)$. If $X = \bot_A$ or $Y = \bot_B$ then $(X,Y)$ also belongs to the tensor $\mathcal{T}$. In any case, we obtain $(x,y) \in T$. Thus, $F \times G \subseteq T$ and $(F,G) \in h(T)$. In all, this proves the equation $h(T) = \mathcal{T}$. $\square$

With some effort, the results in this section may be extended to tensor products with an arbitrary (infinite) number of factors (Erné [15]).

## 4. Tensor products of posets and complete lattices

Diverse tensor products for posets and lattices have been introduced in the literature (see, for example, [4, 10, 11, 26, 28, 34, 36]). We shall define tensor products with two additional parameters governing the choice of the joins that should be preserved by the corresponding bimorphisms.

In order to prepare our construction of tensor products, we consider, for an arbitrary partially ordered set (*poset*) $B$, the *down-sets* or *lower sets* of $B$; these are the fixpoints of the *down-closure operator* $\downarrow$ given by

$$\downarrow Y = \{x \in B : x \leq y \text{ for some } y \in Y\} \ (Y \subseteq B).$$

The down-sets form a completely distributive algebraic lattice with unions as joins and intersections as meets, the *down-set lattice* or *Alexandroff completion* $\mathcal{A}B$. While the down-sets are precisely the *unions* of *principal ideals*

$$\downarrow y = \downarrow\{y\} = \{x \in B : x \leq y\} \ (y \in B),$$

the *intersections* of principal ideals are special down-sets, the (*lower*) *cuts*. They form the *Dedekind-MacNeille completion* or *normal completion* $\mathcal{N}B$ [5]. The associated closure operator is the *cut operator* $\Delta \colon \mathcal{P}B \to \mathcal{P}B$, given by

$$\Delta X = \bigcap \{\downarrow y : X \subseteq \downarrow y\} \ (= \downarrow \bigvee X \text{ if } X \text{ has a supremum } \bigvee X)$$

(see e.g. [16]). More generally, let us consider *augmented posets*, i.e. pairs $B_{\mathcal{Y}} = (B, \mathcal{Y})$ where $B$ is a poset and $\mathcal{Y}$ is a collection of subsets of $B$. We call a subset $I$ of $B$ a $\mathcal{Y}$-*ideal* if $\Delta Y \subseteq I$ for all $Y \in \dot{\mathcal{Y}}$ with $Y \subseteq I$, where

$$\dot{\mathcal{Y}} = \mathcal{Y} \cup \{\{x\} : x \in B\}.$$

The $\mathcal{Y}$-ideals are the fixpoints of the preclosure operator $\Delta_{\mathcal{Y}}$ defined by

$$\Delta_{\mathcal{Y}} X = \bigcup \{\Delta Y : Y \subseteq X, Y \in \dot{\mathcal{Y}}\}.$$

They form a *standard completion* $\mathcal{I}_{\mathcal{Y}} B$, i.e. a closure system between $\mathcal{N}B$ and $\mathcal{A}B$. The least $\mathcal{Y}$-ideal is $\overline{\emptyset}$, the *bottom* of $B$; it is empty or a singleton.

Applying these definitions to any *subset selection* $\mathcal{Z}$, assigning to each poset $B$ a collection $\mathcal{Y} = \mathcal{Z}B$ of subsets, one obtains the $\mathcal{Z}$-ideals as introduced in [8] (see also [14, 16, 18]). A poset $B$ in which each $Y \in \mathcal{Y}$ (resp. $\mathcal{Z}B$) has a join is said to be $\mathcal{Y}$-(resp. $\mathcal{Z}$-) *join-complete*. In such posets, a $\mathcal{Y}$-ideal is just a down-set closed under forming joins of sets in $\mathcal{Y}$ (so-called $\mathcal{Y}$-*joins*).

| $\mathcal{Z}$ | members | $\mathcal{Z}$-join-complete | $\mathcal{Z}$-ideal |
|---|---|---|---|
| $\mathcal{A}$ | arbitrary down-sets | complete lattice | cut [5, 16] |
| $\mathcal{B}$ | binary subsets | join-semilattice | 3-ideal [17] |
| $\mathcal{C}$ | nonempty chains | chain-complete | chain-closed set |
| $\mathcal{D}$ | directed subsets | up-complete (dcpo) | Scott-closed set [7] |
| $\mathcal{E}$ | one-element subsets | arbitrary poset | arbitrary subset |
| $\mathcal{F}$ | finite subsets | join-semilattice with 0 | Frink ideal [22] |

The complete lattices are those posets in which cuts are principal ideals; sometimes, we regard them as augmented posets $(C, \mathcal{P}C)$. In dcpos, the Scott-closed sets are the usual ones [25]; in the absence of directed joins, they are the closed sets in the *weak Scott topology* $\sigma_2 B$ [7]. The Axiom of Choice resp. Zorn's Lemma guarantees that $\mathcal{C}$-join-completeness is equivalent to $\mathcal{D}$-join-completeness, and that the chain-closed sets are the Scott-closed ones.

Now, consider two posets $A$ and $B$ and their direct product $A \times B$ (with the componentwise order). Its down-sets are the *lower relations*, i.e., those relations $R \subseteq A \times B$ which coincide with their down-closure

$$\downarrow R = \{(a,b) \in A \times B : a \leq x \text{ and } b \leq y \text{ for some } (x,y) \in R\}.$$

It is a convenient custom to write $x\,R\,y$ for $(x,y) \in R$, and to put

$$xR = \{y \in B : x\,R\,y\} \text{ for } x \in A, \ XR = \bigcup\{xR : x \in X\} \text{ for } X \subseteq A,$$
$$Ry = \{x \in A : x\,R\,y\} \text{ for } y \in B, \ RY = \bigcup\{Ry : y \in Y\} \text{ for } Y \subseteq B.$$

**Lemma 4.1.** *For any poset $B$, the Alexandroff completion $\mathcal{Q}B = \mathcal{A}(B \times B)$ is a subquantale of the relation quantale $\mathcal{P}(B \times B)$ (but not a quantic quotient unless the order of $B$ is equality). For $R, S, T \in \mathcal{Q}B$, the residuals are given by*

$$R \to T = \{(a,b) \in B \times B : R(\downarrow a) \times \{b\} \subseteq T\},$$
$$T \leftarrow S = \{(a,b) \in B \times B : \{a\} \times (\downarrow b)S \subseteq T\}.$$

*Proof.* Like any down-set lattice, $\mathcal{A}(B \times B)$ is closed under arbitrary unions; these are joins in $\mathcal{P}(B \times B)$ as well as in $\mathcal{Q}B$, and products distribute over them. Since $R, S, T$ are down-sets, so are $R \to T$ and $T \leftarrow S$. Furthermore,

$$S \subseteq R \to T \Leftrightarrow (a\,S\,b \Rightarrow R(\downarrow a) \times \{b\} \subseteq T)$$
$$\Leftrightarrow \forall c, d \ (d \leq a,\ a\,S\,b,\ c\,R\,d \Rightarrow c\,T\,b) \qquad (S = \downarrow S)$$
$$\Leftrightarrow \forall c \ (a\,S\,b,\ c\,R\,a \Rightarrow c\,T\,b) \Leftrightarrow R \cdot S \subseteq T.$$

An analogous argument yields $R \subseteq T \leftarrow S \Leftrightarrow R \cdot S \subseteq T$. □

Notice that the dual relations $R^{op} = \{(y,x) : x\,R\,y\}$ fulfil the equations

$$(R \cdot S)^{op} = S^{op} \cdot R^{op} \text{ and } (R \to T)^{op} = T^{op} \leftarrow R^{op}.$$

Generalizing the tensor product from [10], we now define the tensor product of augmented posets $A_\mathcal{X}$ and $B_\mathcal{Y}$, or the $\mathcal{XY}$-*tensor product* of $A$ and $B$, to be the closure system of all $\mathcal{XY}$-*tensors*,

$$A_\mathcal{X} \otimes B_\mathcal{Y} = \{T = \downarrow T \subseteq A \times B : X \in \mathcal{X}, Y \in \mathcal{Y}, X \times Y \subseteq T \Rightarrow \Delta X \times \Delta Y \subseteq T\}.$$

If $A$ and $B$ happen to be complete lattices or at least all members of $\mathcal{X}$ and $\mathcal{Y}$ possess joins then $T$ is an $\mathcal{XY}$-tensor iff it is a down-set in $A \times B$ so that

$$X \in \mathcal{X},\ Y \in \mathcal{Y} \text{ and } X \times Y \subseteq T \text{ imply } (\bigvee X, \bigvee Y) \in T.$$

Note: $\bigvee(X \times Y) = (\bigvee X, \bigvee Y)$ holds only if $X \neq \emptyset \neq Y$ or $\bigvee X = 0_A, \bigvee Y = 0_B$.

If $A_\mathcal{X}$ and $B_\mathcal{Y}$ are augmented posets then for $T$ to be an $\mathcal{XY}$-tensor it is necessary and sufficient that the following two implications are fulfilled:

$$X \in \dot{\mathcal{X}} \text{ and } X \times \{y\} \subseteq T \text{ imply } \Delta X \times \{y\} \subseteq T,$$
$$Y \in \dot{\mathcal{Y}} \text{ and } \{x\} \times Y \subseteq T \text{ imply } \{x\} \times \Delta Y \subseteq T,$$

which simply means that each $Ty$ is an $\mathcal{X}$-ideal and each $xT$ is a $\mathcal{Y}$-ideal. Hence, the $\mathcal{XY}$-tensor products may be regarded as special instances of tensor products for closure spaces: $A_\mathcal{X} \otimes B_\mathcal{Y}$ is just the tensor product of the closure spaces $(A, \mathcal{I}_\mathcal{X} A)$ and $(B, \mathcal{I}_\mathcal{Y} B)$, and so all facts concerning tensor products of closure spaces apply to that situation. In particular, the tensor product of augmented posets may be reduced to that of complete lattices, by virtue of the following special instance of Theorem 3.2:

**Corollary 4.1.** $A_\mathcal{X} \otimes B_\mathcal{Y} \simeq \mathcal{I}_\mathcal{X} A \otimes \mathcal{I}_\mathcal{Y} B \simeq \mathcal{I}_\mathcal{Y} B \otimes \mathcal{I}_\mathcal{X} A \simeq B_\mathcal{Y} \otimes A_\mathcal{X}.$

**Examples.** Let $A_\mathcal{X}$ and $B_\mathcal{Y}$ be augmented posets.

(1) For $\mathcal{X} = \mathcal{Y} = \emptyset$ (or at least $\mathcal{X} \subseteq \mathcal{E}A$ and $\mathcal{Y} \subseteq \mathcal{E}B$), we have the entire down-set lattice

$$A_\mathcal{X} \otimes B_\mathcal{Y} = A \textcircled{\downarrow} B = \mathcal{A}(A \times B) = \{T \subseteq A \times B : T = {\downarrow}T\}.$$

(2) Taking, at the other extreme, for $\mathcal{X}$ and $\mathcal{Y}$ the power sets $\mathcal{P}A$ and $\mathcal{P}B$, respectively, and identifying $A_{\mathcal{P}A}$ with $A$ etc., we obtain the tensor product

$$A_\mathcal{X} \otimes B_\mathcal{Y} = A \otimes B = \{T \subseteq A \times B : X \times Y \subseteq T \Rightarrow \Delta X \times \Delta Y \subseteq T\}$$

considered in [10]; if $A$ and $B$ are complete lattices, it agrees with the tensor product studied earlier by Shmuely [34] and others; but if $A$ or $B$ fails to be complete, the tensor product $A \otimes B$ will be larger than the tensor product $A_\ell \otimes_r B$ mentioned in the introduction, which consists of all down-sets $T$ in $A \times B$ such that each of the sets $xT$ and $Ty$ is a principal ideal. In [10], it is shown that $A \otimes B$ has much better properties than the tensor product $A_\ell \otimes_r B \simeq Gal(A, B)$ investigated by Nelson [28].

(3) If $\mathcal{X}$ is empty and $\mathcal{Y}$ is the power set $\mathcal{P}B$ then

$$A_\mathcal{X} \otimes B_\mathcal{Y} = A_\emptyset \otimes B = \{T \subseteq A \times B : T = {\downarrow}T, Y \subseteq xT \Rightarrow \Delta Y \subseteq xT\}$$

consists of all down-sets $T$ in $A \times B$ whose slices $xT$ are cuts. Similarly, for $\mathcal{X} = \mathcal{P}A$ and $\mathcal{Y} = \emptyset$ one obtains the system $A \otimes B_\emptyset$ of all down-sets $T$ in $A \times B$ such that each $Ty$ is a cut. If $A$ and $B$ are complete lattices, $A_\emptyset \otimes B$ is the system of all right tensors, $A \otimes_r B \simeq Ant(A, B)$, and $A \otimes B_\emptyset$ is the system of all left tensors, $A_\ell \otimes B \simeq Ant(B, A)$ (see the introduction).

(4) A variant of the previous example is obtained for $\mathcal{X} = \{\emptyset\}$ and $\mathcal{Y} = \mathcal{P}B$. In case $B$ is complete, the members of $A_\mathcal{X} \otimes B_\mathcal{Y}$ now correspond to the antitone maps sending the least element of $A$ to the greatest element of $B$.

(5) Taking $\mathcal{X} = \{X \subseteq A : card\, X < \alpha\}$ and $\mathcal{Y} = \{Y \subseteq B : card\, Y < \alpha\}$, one obtains some (but not all) kinds of tensor products studied by Bandelt [4].

## 5. Morphisms and bimorphisms for posets

The introductory remarks on the isomorphisms

$$f \mapsto T_f = \{(x,y) : f(x) \geq y\} \text{ and } T \mapsto f_T \text{ with } f_T(x) = \max xT$$

may be generalized as follows. Let $A_{\mathcal{X}}$ be an augmented poset, $B$ a poset, and define $A_{\mathcal{X}} \otimes B$ to be the set of all $T \subseteq A \times B$ such that for all $x \in A$, the set $xT$ is a principal ideal, and for all $y \in B$, the set $Ty$ is an $\mathcal{X}$-ideal. For $Ant_{\mathcal{X}}(A,B)$ take the set of all maps $f \colon A \to B$ such that for each $y \in B$, the preimage $f^{-1}[\uparrow y]$ is an $\mathcal{X}$-ideal. Thus, $Ant_{\mathcal{X}}(A,B)$ consists of all $\mathcal{X}$-ideal continuous maps from $A$ to the order-dual of $B$, where a map $f$ between posets is $\mathcal{X}$-*ideal continuous* if all preimages of principal ideals are $\mathcal{X}$-ideals. For complete lattices $A, B$, one obtains the complete lattices $A_{\mathcal{X}} \otimes B$ and

$$Ant_{\mathcal{X}}(A,B) = \{f \in Ant(A,B) : f(\bigvee X) = \bigwedge f[X] \text{ for all } X \in \mathcal{X}\}.$$

**Proposition 5.1.** *For any two posets $A$ and $B$, the mutually inverse isomorphisms $f \mapsto T_f$ and $T \mapsto f_T$ between the posets $Ant(A,B)$ and $A \otimes_r B$ induce mutually inverse isomorphisms between $Ant_{\mathcal{X}}(A,B)$ and $A_{\mathcal{X}} \otimes B$.*

*Proof.* By definition, we have for $f \in Ant(A,B)$ and $T \in A \otimes_r B$:

$$T_f \subseteq T \Leftrightarrow \forall x \in A \ (x \, T \, f(x)) \Leftrightarrow \forall x \in A \ (f(x) \leq f_T(x)) \Leftrightarrow f \leq f_T.$$

For $f \in Ant_{\mathcal{X}}(A,B)$ and each $x \in A$, the set $xT_f$ is a principal ideal $\downarrow f(x)$, and for each $y \in B$, the set $T_f \, y = f^{-1}[\uparrow y]$ is an $\mathcal{X}$-ideal. Thus, we get $T_f \in A_{\mathcal{X}} \otimes B$ and $f(x) = \max xT_f$.

Conversely, let $T \in A_{\mathcal{X}} \otimes B$. Then $f_T(x) = \max xT$ defines an antitone map by the down-set property of $T$. For $x \in A$ and $y \in B$, we have the equivalences

$$x \in f_T^{-1}[\uparrow y] \Leftrightarrow f_T(x) \geq y \Leftrightarrow x \, T \, y.$$

Hence, $f_T^{-1}[\uparrow y] = Ty$ is an $\mathcal{X}$-ideal, proving $f_T \in Ant_{\mathcal{X}}(A,B)$.

The equivalence $x \, T \, y \Leftrightarrow f_T(x) \geq y$ shows that the assignment $T \mapsto f_T$ induces a bijection between $A_{\mathcal{X}} \otimes B$ and $Ant_{\mathcal{X}}(A,B)$ with inverse $f \mapsto T_f$. Since $f \leq g$ is tantamount to $T_f \subseteq T_g$, we actually have an isomorphism. □

Let $A_{\mathcal{X}}$ and $B_{\mathcal{Y}}$ be augmented posets. Being a closure system, $A_{\mathcal{X}} \otimes B_{\mathcal{Y}}$ is a complete lattice with greatest element $A \times B$ and least element $\overline{\emptyset}$, which is one of the following sets:

$$
\begin{array}{lll}
\text{if} & \emptyset \notin \mathcal{X}, \emptyset \notin \mathcal{Y} & \text{then} \quad \overline{\emptyset} = \emptyset \\
\text{if} & \emptyset \in \mathcal{X}, \emptyset \notin \mathcal{Y} & \text{then} \quad \overline{\emptyset} = \Delta\emptyset \times B \\
\text{if} & \emptyset \notin \mathcal{X}, \emptyset \in \mathcal{Y} & \text{then} \quad \overline{\emptyset} = A \times \Delta\emptyset \\
\text{if} & \emptyset \in \mathcal{X}, \emptyset \in \mathcal{Y} & \text{then} \quad \overline{\emptyset} = A \times \Delta\emptyset \cup \Delta\emptyset \times B.
\end{array}
$$

Here, $\Delta\emptyset$ is empty or a singleton. Down-sets of the form

$$x \otimes y = \downarrow(x,y) \cup \overline{\emptyset} = \bigcap \{T \in A_{\mathcal{X}} \otimes B_{\mathcal{Y}} : x \, T \, y\}$$

are always members of $A_{\mathcal{X}} \otimes B_{\mathcal{Y}}$; we call them *pure tensors* or *principal tensors*. They form a join-dense subset of $A_{\mathcal{X}} \otimes B_{\mathcal{Y}}$, because every $\mathcal{X}\mathcal{Y}$-tensor $T$ is the union (hence the join) of all pure tensors $x \otimes y$ with $x \, T \, y$.

As expected, $\mathcal{X}\mathcal{Y}$-tensor products may be characterized by a universal property for bimorphisms (cf. [2, 10, 11]). Notice that a map from a poset $A$ into a complete lattice $C$ is $\mathcal{X}$-ideal continuous iff it is continuous as a map between the closure spaces $(A, \mathcal{I}_\mathcal{X} A)$ and $(C, \mathcal{N}C) = (C, \mathcal{I}_{\mathcal{P}C} C)$. Given a further augmented poset $B_\mathcal{Y}$, we call a map $f\colon A \times B \to C$ an $\mathcal{X}\mathcal{Y}$-*bimorphism* if each $f_y$ is $\mathcal{X}$-ideal continuous and each ${}_x f$ is $\mathcal{Y}$-ideal continuous — in other words, if $f$ is separately continuous as a map between the associated closure spaces. In the case of complete lattices $A$ and $B$, $f$ is an $\mathcal{X}\mathcal{Y}$-bimorphism iff it is isotone and satisfies

$${}_x f(\bigvee Y) = \bigvee {}_x f[Y] \text{ for all } x \in A,\ Y \in \mathcal{Y},$$

$$f_y(\bigvee X) = \bigvee f_y[X] \text{ for all } y \in B,\ X \in \mathcal{X}.$$

Now, as a special instance of the corresponding Theorem 3.1, we have:

**Corollary 5.1.** *For arbitrary augmented posets $A_\mathcal{X}$ and $B_\mathcal{Y}$, the map*

$$\otimes\colon A \times B \to A_\mathcal{X} \otimes B_\mathcal{Y},\ (x,y) \mapsto x \otimes y$$

*is universal among all $\mathcal{X}\mathcal{Y}$-bimorphisms into complete lattices. Thus, $\otimes$ is an $\mathcal{X}\mathcal{Y}$-bimorphism, and for any $\mathcal{X}\mathcal{Y}$-bimorphism $f$ with domain $A \times B$, there is a unique join-preserving map*

$$f^\vee\colon A_\mathcal{X} \otimes B_\mathcal{Y} \to C$$

*with $f(x,y) = f^\vee(x \otimes y)$.*

**Examples.** For simplicity, let $A$ and $B$ be complete lattices.

(1) The down-set lattice $\mathcal{A}(A \times B) = A_\emptyset \otimes B_\emptyset$ is the tensor product for isotone bimorphisms.

(2) $A \otimes B = A_{\mathcal{P}A} \otimes B_{\mathcal{P}B} \simeq Gal(A,B)$ is the known tensor product for bimorphisms that preserve arbitrary joins in both arguments.

(3) $A_\emptyset \otimes B \simeq Ant(A,B)$ is the tensor product for bimorphisms that are isotone in the first argument and preserve joins in the second argument.

(4) For $\mathbf{1} = \{\emptyset\}$, $A_\mathbf{1} \otimes B$ is slightly different from $A_\emptyset \otimes B$, being the tensor product for bimorphisms that are isotone, preserve bottom elements in the first argument, and joins in the second. Here, $A_\mathbf{1} \otimes B$ is isomorphic to $Ant_1(A,B)$, the lattice of all antitone maps $f\colon A \to B$ sending $0_A$ to $1_B$.

(5) If $\mathcal{X}$ and $\mathcal{Y}$ consist of all finite subsets then $A_\mathcal{X} \otimes B_\mathcal{Y}$ is the tensor product for bimorphisms preserving finite joins in both arguments. The same universal property holds for bimorphisms preserving joins of cardinality less than an arbitrary but fixed cardinal $\alpha$, if we take for $\mathcal{X}$ and $\mathcal{Y}$ the collections of all subsets with less than $\alpha$ elements.

(6) If $\mathcal{X}$ and $\mathcal{Y}$ consist of all directed subsets then $A_\mathcal{X} \otimes B_\mathcal{Y}$ is the tensor product for Scott-continuous bimorphisms (cf. Bandelt [3]).

## 6. Truncated tensor products

For subsets $X$ of a closure space $A = (E, \mathcal{C})$ with $\perp_A = \overline{\emptyset}$, put $\check{X} = X \smallsetminus \perp_A$ and denote by $\check{A}$ the *unbounded coreflection* of $A$, the subspace induced on $\check{E}$. We shall now present a more economic and effective form of tensor products, by cutting off from each tensor the obligatory part $\overline{\emptyset} = (\perp_A \times \mathcal{U}B) \cup (\mathcal{U}A \times \perp_B)$.

**Proposition 6.1.** *If $A$ and $B$ are closure spaces then $A \otimes B$ is isomorphic to $\check{A} \otimes \check{B}$, by the mutually inverse assignments $T \mapsto T \smallsetminus \overline{\emptyset}$ and $T' \mapsto T' \cup \overline{\emptyset}$.*

*Proof.* $T \in A \otimes B$ implies $T \smallsetminus \overline{\emptyset} \in \check{A} \otimes \check{B}$ and $(T \smallsetminus \overline{\emptyset}) \cup \overline{\emptyset} = T$: $X \times Y \subseteq T \smallsetminus \overline{\emptyset}$ entails $\overline{X} \times \overline{Y} \subseteq T$, hence $(\overline{X} \smallsetminus \perp_A) \times (\overline{Y} \smallsetminus \perp_B) = (\overline{X} \times \overline{Y}) \smallsetminus \overline{\emptyset} \subseteq T \smallsetminus \overline{\emptyset}$.

Conversely, $T' \in \check{A} \otimes \check{B}$ implies $T = T' \cup \overline{\emptyset} \in A \otimes B$ and $T' = T \smallsetminus \overline{\emptyset}$: from $X \times Y \subseteq T$ it follows that $(X \smallsetminus \perp_A) \times (Y \smallsetminus \perp_B) = (X \times Y) \smallsetminus \overline{\emptyset} \subseteq T \smallsetminus \overline{\emptyset} = T'$, hence $(\overline{X} \times \overline{Y}) \smallsetminus \overline{\emptyset} = (\overline{X} \smallsetminus \perp_A) \times (\overline{Y} \smallsetminus \perp_B) \subseteq T'$ and so $\overline{X} \times \overline{Y} \subseteq T' \cup \overline{\emptyset} = T$. $\square$

Given an augmented poset $A_\mathcal{X} = (A, \mathcal{X})$ with least $\mathcal{X}$-ideal $\perp_A$, put

$$\check{\mathcal{X}} = \{\check{X} = X \smallsetminus \perp_A : X \in \dot{\mathcal{X}},\ \check{X} \neq \emptyset \text{ or } \perp_A = \emptyset\}.$$

Note that the cut operator $\check{\Delta}$ of the subposet $\check{A}$ satisfies

$$\check{\Delta}\check{X} \supseteq \bigcap \{\downarrow y \smallsetminus \perp_A : \check{X} \subseteq \downarrow y\} = \bigcap \{\downarrow y \smallsetminus \perp_A : X \subseteq \downarrow y\} = \Delta X \smallsetminus \perp_A$$

and equality holds if $\check{X} \in \check{\mathcal{X}}$. Using that fact, one easily shows:

**Proposition 6.2.** *The ideal lattices $\mathcal{I}_\mathcal{X} A$ and $\mathcal{I}_{\check{\mathcal{X}}} \check{A}$ are isomorphic, by virtue of the mutually inverse assignments $I \mapsto I \smallsetminus \overline{\emptyset}$ and $I' \mapsto I' \cup \overline{\emptyset}$.*

For augmented posets $A_\mathcal{X}$ and $B_\mathcal{Y}$, the *truncated tensor product*

$$A_\mathcal{X} \check{\otimes} B_\mathcal{Y} = \check{A}_{\check{\mathcal{X}}} \otimes \check{B}_{\check{\mathcal{Y}}}$$

consists of all subsets $T$ of the product poset $\check{A} \times \check{B}$ such that each $xT$ is a $\check{\mathcal{Y}}$-ideal and each $Ty$ is an $\check{\mathcal{X}}$-ideal. Thus, $A_\mathcal{X} \check{\otimes} B_\mathcal{Y}$ is a closure system with least element $\emptyset$. The least element of $A_\mathcal{X} \check{\otimes} B_\mathcal{Y}$ containing $(x, y)$ is the principal ideal $\downarrow(x, y)$ in $\check{A} \times \check{B}$. Passing from $A_\mathcal{X}$ and $B_\mathcal{Y}$ to their $\mathcal{X}$- and $\mathcal{Y}$-ideal spaces, respectively, we deduce from Propositions 6.1 and 6.2:

**Corollary 6.1.** *For augmented posets $A_\mathcal{X}, B_\mathcal{Y}$, the tensor product $A_\mathcal{X} \otimes B_\mathcal{Y}$ is isomorphic to $A_\mathcal{X} \check{\otimes} B_\mathcal{Y}$ via the relativization map $T \mapsto T \cap (\check{A} \times \check{B}) = T \smallsetminus \overline{\emptyset}$.*

**Corollary 6.2.** *For complete $B$, the map $T \mapsto \check{f}_T$ with $\check{f}_T(x) = \max(xT \cup \overline{\emptyset})$ induces an isomorphism between $A_\mathcal{X} \check{\otimes} B = A_\mathcal{X} \check{\otimes} B_{\mathcal{P}B}$ and $Ant_\mathcal{X}(A, B)$. The inverse isomorphism is given by $f \mapsto \check{T}_f = \{(x, y) \in \check{A} \times \check{B} : f(x) \geq y\}$.*

**Corollary 6.3.** *For posets $A, B$, the principal ideal map $(x, y) \mapsto \downarrow(x, y)$ from $A \times B$ to $A \check{\otimes} B$ is a universal bimorphism into complete lattices.*

## 7. Pseudocomplementation and polarization

The *polar* $x^\perp$ of an element $x$ in a poset $B$ consists of all $y \in B$ such that $\downarrow x \cap \downarrow y = \Delta\emptyset$. If the polar $x^\perp$ has a greatest element, this is denoted by $x^*$ and called the *pseudocomplement* of $x$; and the whole poset is said to be *pseudocomplemented* if each of its elements has a pseudocomplement.

Recall that $A_\ell \otimes_r B$ denotes the set of all tensors $T \subseteq A \times B$ so that each $xT$ and each $Ty$ is a principal ideal. For bounded posets, pseudocomplementation is preserved by these tensor products. Write $A \otimes b$ for $(A \times \downarrow b) \cup \overline{\emptyset}$, the least tensor containing $A \times \downarrow b$, and $a \otimes B$ for $(\downarrow a \times B) \cup \overline{\emptyset}$.

**Proposition 7.1.** *Let $A$ and $B$ be non-singleton bounded posets.*

(1) *For $b \in B$, a set $S$ is the pseudocomplement of $A \otimes b$ in $A_\ell \otimes_r B$ iff $b$ has a pseudocomplement $b^*$ in $B$ and $S = A \otimes b^*$.*

(2) $A_\ell \otimes_r B$ *is pseudocomplemented iff $A$ and $B$ are pseudocomplemented.*

*Proof.* (1) If $b^*$ is the pseudocomplement of $b$ then $A \otimes b^*$ is a tensor with $A \otimes b \cap A \otimes b^* = A \times (\downarrow b \cap \downarrow b^*) \cup \overline{\emptyset} = \overline{\emptyset}$; for $T \in A_\ell \otimes_r B$ with $A \otimes b \cap T = \overline{\emptyset}$ and any $(x, y) \in T$, it follows that $\overline{\emptyset} = A \otimes b \cap x \otimes y = \downarrow x \times (\downarrow y \cap \downarrow b) \cup \overline{\emptyset}$, hence $x = 0_A$ or $y \in b^\perp$; in any case, $(x, y) \in (\downarrow 0_A \times B) \cup (A \times \downarrow b^*) \subseteq A \otimes b^*$; thus, $T \subseteq A \otimes b^*$. Consequently, $A \otimes b^*$ is the pseudocomplement of $A \otimes b$.

Conversely, if $S$ is the pseudocomplement of $A \otimes b$ in $A_\ell \otimes_r B$ then for any $a \in \breve{A} \neq \emptyset$ (by hypothesis, $A$ is not a singleton), the set $aS$ is a principal ideal, so it suffices to prove $aS = b^\perp$ in order to show that $b$ has a pseudocomplement $b^*$ (the greatest element of $aS = b^\perp$): $aSy$ implies $\downarrow a \times (\downarrow b \cap \downarrow y) \subseteq a \otimes b \cap a \otimes y = \overline{\emptyset}$, hence $\downarrow b \cap \downarrow y = \Delta\emptyset$, i.e. $b \perp y$. Conversely, if the latter holds then $A \otimes b \cap a \otimes y = \downarrow a \times (\downarrow b \cap \downarrow y) \cup \overline{\emptyset} = \overline{\emptyset}$, hence $(a, y) \in a \otimes y \subseteq S$, $y \in aS$.

(2) If $A$ and $B$ are pseudocomplemented then, by (1), each principal tensor $a \otimes b$ has the pseudocomplement $S = A \otimes b^* \cup a^* \otimes B$. Indeed, $S$ belongs to $A_\ell \otimes_r B$, as $xS \in \{\downarrow b^*, \downarrow 1_B\}$ for all $x$ in $A$, and $Sy \in \{\downarrow a^*, \downarrow 1_A\}$ for all $y$ in $B$. Furthermore, $a \otimes b \cap S = a \otimes 0_B \cup 0_A \otimes b = \overline{\emptyset}$, and if a $T \in A_\ell \otimes_r B$ satisfies $a \otimes b \cap T = \overline{\emptyset}$ then for each $(x, y) \in T$ we get $\overline{\emptyset} = x \otimes y \cap a \otimes b$, hence $x \leq a^*$ or $y \leq b^*$, in any case $(x, y) \in S$. Since the set of all principal tensors is join-dense in $A_\ell \otimes_r B$, all elements of $A_\ell \otimes_r B$ have pseudocomplements. Indeed, if $J$ is a join-dense subset of a complete lattice $L$ and all $x \in J$ have pseudocomplements $x^*$ then each $t \in L$ has the pseudocomplement $t^* = \bigwedge \{x^* : x \in J \cap \downarrow t\}$.

Conversely, if $A_\ell \otimes_r B$ is pseudocomplemented then by (1), so are $A$ and $B$. □

Let us call a closure space $A$ *polarized* if for each $x \in \mathcal{U}A$, the *polar*

$$x^\perp = \{y \in \mathcal{U}A : \overline{x} \cap \overline{y} = \overline{\emptyset}\}$$

is closed. In particular, a complete lattice is pseudocomplemented iff it is polarized as a closure space. These specific notions of polars fit into the general framework of polarities proposed by Birkhoff [5] (cf. [16, 17, 18]).

**Theorem 7.1.** *For closure spaces $A, B$, the following are equivalent:*

(a) *$A$ and $B$ are polarized, or one of them has only one closed set.*
(b) *The complete lattices $\mathcal{C}A$ and $\mathcal{C}B$ are pseudocomplemented, or one of them is a singleton.*
(c) *The lattices $\check{A} \otimes \check{B} \simeq A \otimes B \simeq \mathcal{C}A \otimes \mathcal{C}B$ are pseudocomplemented.*

*Proof.* (a) $\Rightarrow$ (b): For $F \in \mathcal{C}A$, the polar $F^\perp = \bigcap \{y^\perp : y \in F\}$ is closed, and for $G \in \mathcal{C}B$, the equation $F \cap G = \overline{\emptyset}$ is equivalent to $G \subseteq F^\perp$.

(b) $\Rightarrow$ (a): Each polar $x^\perp$ is contained in the pseudocomplement $K$ of $\overline{x}$ in $\mathcal{C}A$; hence $\overline{x^\perp} \subseteq K$. Thus, if $y \in \overline{x^\perp}$ then $y \in x^\perp$, which is therefore closed. For (b) $\Leftrightarrow$ (c), invoke Theorem 3.2 and Propositions 6.1 and 7.1. $\square$

A poset $B$ is said to be $\mathcal{Y}$-*distributive at the bottom* or $\perp$-$\mathcal{Y}$-*distributive* provided for all $x \in B$ and $Y \in \mathcal{Y}$, $x \perp Y$ implies $x \perp \Delta Y$, where $x \perp Y$ means that $\downarrow x \cap \downarrow Y$ is contained in $\overline{\emptyset}$, the least $\mathcal{Y}$-ideal. Writing $x \perp y$ for $x \perp \{y\}$, we see that $\perp$-$\mathcal{Y}$-distributivity simply requires that all polars $x^\perp = \{y : x \perp y\}$ are $\mathcal{Y}$-ideals. An equivalent postulate is that

$$a \in \Delta Y \cap \check{B} \text{ implies } \downarrow a \cap \downarrow Y \cap \check{B} \neq \emptyset \text{ for all } Y \in \mathcal{Y}.$$

The term $\perp$-$\mathcal{Y}$-*distributive* is motivated by the fact that a ($\mathcal{Y}$-join) complete lattice $B$ is $\perp$-$\mathcal{Y}$-distributive iff $x \wedge y = 0_B$ for all $y \in Y$ implies $x \wedge \bigvee Y = 0_B$. For more material on distributivity at the bottom see Erné and Joshi [17].

**Corollary 7.1.** *A poset $B$ is $\perp$-$\mathcal{Y}$-distributive iff $\mathcal{I}_\mathcal{Y} B$ is pseudocomplemented.*

*Proof.* By definition, $B$ is $\perp$-$\mathcal{Y}$-distributive iff the ideal space $(B, \mathcal{I}_\mathcal{Y} B)$ is polarized, which by Theorem 7.1 says that $\mathcal{I}_\mathcal{Y} B$ is pseudocomplemented. $\square$

**Examples.** Let $B$ be a poset.

(1) If $\mathcal{Y} \subseteq \mathcal{E}B \cup \{\emptyset\}$ then $B$ is trivially $\perp$-$\mathcal{Y}$-distributive.

(2) At the other extreme, for the power set $\mathcal{Y} = \mathcal{P}B$, the $\perp$-$\mathcal{Y}$-distributive posets are just the *completely $\perp$-distributive* ones, in which all polars are cuts. Pseudocomplemented posets are completely $\perp$-distributive; for complete lattices, the converse is also true. Clearly, frames are pseudocomplemented.

(3) If $\mathcal{Y}$ is the collection of all *finite* subsets of a join-semilattice with 0, then $\perp$-$\mathcal{Y}$-distributivity amounts to the usual notion of 0-(*semi-*)*distributivity*.

(4) If $\mathcal{Y}$ consists of all *directed subsets*, one obtains complete lattices that are *meet-continuous at* 0. That property together with atomicity (requiring an atom below any non-zero element) has the consequence that the least element is a meet of $\bigwedge$-irreducible elements, a fact of great importance for universal algebra and spectral theory. For an easy proof, observe that any polar $a^\perp$ is closed under directed joins, hence has maximal elements by Zorn's Lemma, and these are $\bigwedge$-irreducible if $a$ is an atom, because every greater element has to lie above $a$. Note that a complete lattice is pseudo-complemented iff it is 0-distributive and meet-continuous at 0.

## 8. Tensor quantales

A first observation concerning the structure of truncated tensor products is

**Lemma 8.1.** *For any poset B, the truncated tensor product*

$$\mathcal{Q}\check{B} = \check{B} \mathbin{\check{\mathbb{Q}}} \check{B} = \mathcal{A}(\check{B}\times\check{B}) = B_\emptyset \mathbin{\check{\otimes}} B_\emptyset$$

*is a quantale with union as join, intersection as meet, and the relation product as multiplication. It has a unit only if the order relation on $\check{B}$ is equality.*

For the last claim, modify the proof of Proposition 9.1 suitably.

Let $\mathcal{X}$ and $\mathcal{Y}$ be subsets of $\mathcal{P}B$. We define on each truncated tensor product $B_\mathcal{X} \mathbin{\check{\otimes}} B_\mathcal{Y}$, which is a meet-closed subset of $\mathcal{Q}\check{B}$, a multiplication $\odot$ by taking for $R \odot S$ the least member of $B_\mathcal{X} \mathbin{\check{\otimes}} B_\mathcal{Y}$ containing the relation product $R \cdot S$. This gives an interesting and flexible multiplicative structure, in contrast to the corresponding multiplication on the "full" tensor products $B_\mathcal{X} \otimes B_\mathcal{Y}$, which becomes trivial whenever $B$ has a least element $0_B$: for $R, S \in B_\mathcal{X} \otimes B_\mathcal{Y}$, the inclusions $B \times \{0_B\} \subseteq R$ and $\{0_B\} \times B \subseteq S$ force $R \cdot S$ and $R \odot S$ to be the largest tensor $B \times B$, no matter how small $R$ and $S$ are.

The crucial question is now: under what circumstances is the truncated tensor product $B_\mathcal{X} \mathbin{\check{\otimes}} B_\mathcal{Y}$ a quantic quotient of $\mathcal{Q}\check{B}$, hence a quantale in its own right? As we know from Section 2, this is the case precisely when we can represent $B_\mathcal{X} \mathbin{\check{\otimes}} B_\mathcal{Y}$ as the fixpoint set of a (pre)nucleus.

In the more general situation of truncated tensor products of augmented posets $A_\mathcal{X}$ and $B_\mathcal{Y}$, define a preclosure operation $t$ on $\check{A} \mathbin{\check{\mathbb{Q}}} \check{B} = \mathcal{A}(\check{A}\times\check{B})$ by

$$t(R) = t_{\mathcal{X}\mathcal{Y}}(R) = \bigcup \{\check{\Delta}X \times \check{\Delta}Y : X \in \check{\mathcal{X}}, Y \in \check{\mathcal{Y}}, X \times Y \subseteq R\}.$$

The cut operators $\check{\Delta}$ refer to the truncated posets $\check{A}$ and $\check{B}$. By definition, the truncated tensor product $A_\mathcal{X} \mathbin{\check{\otimes}} B_\mathcal{Y}$ is the fixpoint set of $t$. This preclosure operation is not idempotent (except in some extremal cases), so that the least closure operation $\overline{t}$ above $t$ has to be built by transfinite iteration, or just by

$$\overline{t}(R) = \overline{t_{\mathcal{X}\mathcal{Y}}}(R) = \bigcap \{T \in A_\mathcal{X} \mathbin{\check{\otimes}} B_\mathcal{Y} : R \subseteq T\}.$$

Specifically, we define new relation products by

$$R \odot S = R_{\mathcal{X}} \odot_{\mathcal{Y}} S = \overline{t}(R \cdot S).$$

Although the usual relation product is certainly associative, it is not clear *a priori* under what conditions the new multiplication $\odot$ is associative as well. We shall see that the appropriate ingredient for positive results is polarization in the case of closure spaces and distributivity at the bottom in the case of posets.

Let us prepare the characterization of those posets for which $B_\mathcal{X} \mathbin{\check{\otimes}} B_\mathcal{Y}$ is a quantale or, more precisely, a quantic quotient of $\mathcal{Q}\check{B}$, by

**Lemma 8.2.** *Let $A_\mathcal{X}, B_\mathcal{Y}, C_\mathcal{Z}$ be augmented posets. If $B$ is $\perp$-$\mathcal{Y}$-distributive then for all $R \in \breve{A} \mathbin{\textcircled{\raisebox{-.1ex}{$\downarrow$}}} \breve{B}$ and $S \in \breve{B} \mathbin{\textcircled{\raisebox{-.1ex}{$\downarrow$}}} \breve{C}$,*
$$t_{\mathcal{X}\mathcal{Y}}(R) \cdot S \cup R \cdot t_{\mathcal{Y}\mathcal{Z}}(S) \subseteq t_{\mathcal{X}\mathcal{Z}}(R \cdot S),$$
$$\overline{t_{\mathcal{X}\mathcal{Y}}}(R) \odot \overline{t_{\mathcal{Y}\mathcal{Z}}}(S) = \overline{t_{\mathcal{X}\mathcal{Z}}}(R \cdot S).$$

*Proof.* For $(x, z) \in t_{\mathcal{X}\mathcal{Y}}(R) \cdot S$, there are $y \in \breve{B}$, $X \in \mathcal{X}$, $Y \in \mathcal{Y}$ with $\breve{X} \neq \emptyset \neq \breve{Y}$, $\breve{X} \times \breve{Y} \subseteq R$, $x \in \breve{\Delta} \breve{X}$, $y \in \breve{\Delta} \breve{Y}$, and $y \, S \, z$. Recall that $\breve{\Delta} \breve{Y} = \Delta Y \smallsetminus \perp_B$, as $\breve{Y} \neq \emptyset$.

By $\perp$-$\mathcal{Y}$-distributivity of $B$, we find a $b \in {\downarrow} y \cap {\downarrow} Y \cap \breve{B}$. It follows that $\breve{X} \times \{b\} \subseteq R$ and $(b, z) \in S$, whence $\breve{X} \times \{z\} \subseteq R \cdot S$, and as $x \in \breve{\Delta} \breve{X}$, one obtains $(x, z) \in t_{\mathcal{X}\mathcal{Z}}(R \cdot S)$. Analogously, one shows $R \cdot t_{\mathcal{Y}\mathcal{Z}}(S) \subseteq t_{\mathcal{X}\mathcal{Z}}(R \cdot S)$.

The equation for the associated closure operations is then obtained by applying Lemma 2.1 to $i = t_{\mathcal{X}\mathcal{Y}}$, $j = t_{\mathcal{Y}\mathcal{Z}}$ and $k = t_{\mathcal{X}\mathcal{Z}}$. $\square$

By similar arguments, one proves $t_{\emptyset \mathcal{Y}}(R) \cdot S = R \cdot S = R \cdot t_{\mathcal{Y} \emptyset}(S)$.

Now, to the main result.

**Theorem 8.1.** *For a poset $B$ and sets $\mathcal{X}, \mathcal{Y} \subseteq \mathcal{P}B$, the following conditions are equivalent:*

(a) $B$ *is* $\perp$-$(\mathcal{X} \cup \mathcal{Y})$-*distributive.*
(b) $t_{\mathcal{X}\mathcal{Y}}$ *is a prenucleus on* $\mathcal{Q}\breve{B}$.
(c) $\overline{t_{\mathcal{X}\mathcal{Y}}}$ *is a nucleus on* $\mathcal{Q}\breve{B}$.
(d) $B_\mathcal{X} \breve{\otimes} B_\mathcal{Y}$ *is a quantic quotient of* $\mathcal{Q}\breve{B}$.
(e) $B_\mathcal{X} \breve{\otimes} B_\mathcal{Y}$ *is a quantale with multiplication* $_\mathcal{X} \odot _\mathcal{Y}$.
(f) $B_\mathcal{X} \breve{\otimes} B_\mathcal{Y}$ *(or $B_\mathcal{X} \otimes B_\mathcal{Y}$) is pseudocomplemented.*
(g) *There exists a (pre)nucleus $j$ on $\mathcal{Q}\breve{B}$ with $t_{\mathcal{X}\mathcal{Y}} \leq j$ and $j(\emptyset) = \emptyset$.*
(h) *There exist $\mathcal{X}$- resp. $\mathcal{Y}$-ideal continuous order embeddings $g$ and $h$ of $B$ in a (pre)quantale $Q$ such that $x \perp y \Leftrightarrow g(x) \cdot h(y) = 0_Q$.*

*Each of these conditions is fulfilled whenever $B$ is pseudocomplemented.*

*Proof.* (a) $\Rightarrow$ (b): Applying Lemma 8.2 to the special case $A = B = C$ and $\mathcal{Y} = \mathcal{Z}$, one obtains for $R, S \subseteq \mathcal{Q}\breve{B}$ first $t_{\mathcal{X}\mathcal{Y}}(R) \cdot S \subseteq t_{\mathcal{X}\mathcal{Y}}(R \cdot S)$, and second $R \cdot t_{\mathcal{X}\mathcal{Y}}(S) \subseteq R \cdot t_{\mathcal{U}\mathcal{Y}}(S) \subseteq t_{\mathcal{X}\mathcal{Y}}(R \cdot S)$ for $\mathcal{U} = \mathcal{X} \cup \mathcal{Y}$.

(a) $\Leftrightarrow$ (f): We may assume that $B$ is not a singleton. Then

$B$ is $\perp$-$(\mathcal{X} \cup \mathcal{Y})$-distributive

$\Leftrightarrow \mathcal{I}_\mathcal{X} B$ and $\mathcal{I}_\mathcal{Y} B$ are pseudocomplemented (Corollary 7.1)

$\Leftrightarrow \mathcal{I}_\mathcal{X} B \otimes \mathcal{I}_\mathcal{Y} B$ is pseudocomplemented (Proposition 7.1)

$\Leftrightarrow B_\mathcal{Y} \otimes B_\mathcal{Y} \simeq B_\mathcal{Y} \breve{\otimes} B_\mathcal{Y}$ is pseudocomplemented (Corollaries 4.1, 6.1).

(b) $\Rightarrow$ (c) $\Leftrightarrow$ (d) $\Leftrightarrow$ (e): Proposition 2.1.

(c) $\Rightarrow$ (g): Take $j = \overline{t_{\mathcal{X}\mathcal{Y}}}$.

(g) $\Rightarrow$ (a): Assume $\breve{B} \neq \emptyset$ (otherwise $B = \perp_B$ is pseudocomplemented). For $Y \in \mathcal{Y}$ and $a \in \Delta Y \cap \breve{B} \subseteq \breve{\Delta} \breve{Y}$, put $R = \breve{B} \times (\breve{B} \cap {\downarrow} \breve{Y})$, $S = (\breve{B} \cap {\downarrow} a) \times \breve{B}$. Then $R, S \in \mathcal{Q}\breve{B}$. If $\breve{Y} \in \breve{\mathcal{Y}}$ then $\breve{B} \times \breve{\Delta} \breve{Y} \subseteq t_{\mathcal{X}\mathcal{Y}}(R) \subseteq j(R)$. Picking a $b \in$ $\breve{B}$, we get $(b,a) \in j(R)$, $(a,b) \in S$, and then $(b,b) \in j(R){\cdot}S \subseteq j(R{\cdot}S)$. Now, by contraposition, the hypothesis $j(\emptyset) = \emptyset$ implies $R{\cdot}S \neq \emptyset$, which provides $x,y,z \in \breve{B}$ such that $x\,R\,y\,S\,z$. Then $y \in \mathord{\downarrow}Y$ (since $(x,y) \in R$) and $y \in \mathord{\downarrow}a$ (since $(y,z) \in S$). Thus, $y \in \breve{B} \cap \mathord{\downarrow}a \cap \mathord{\downarrow}Y$, proving $\bot$-$\mathcal{Y}$-distributivity. It remains to check the case $\breve{Y} = \emptyset$. Then either $\breve{\Delta}\breve{Y} = \emptyset$, and one may proceed as before, or $\Delta\breve{Y} = \{0_{\breve{B}}\}$, and then $B$ is trivially pseudocomplemented, because $\breve{B} = B \smallsetminus \{0_B\}$ has a least element.

(e) $\Rightarrow$ (h): Fix some $b \in \breve{B}$ (the case $\breve{B} = \emptyset$ is obvious). Put $Q = B_{\mathcal{X}} \breve{\otimes} B_{\mathcal{Y}}$ and define $g \colon B \to Q$ by $g(x) = \mathord{\downarrow}(b,x)$ (the down-set in $\breve{B}$) for $x \in \breve{B}$ and $g(x) = \emptyset$ for $x \in \bot_B$. Then $g$ is an order embedding, since for $x,y \in \breve{B}$,

$$x \leq y \Rightarrow g(x) = \mathord{\downarrow}(b,x) \subseteq g(y) = \mathord{\downarrow}(b,y) \Rightarrow (b,x) \leq (b,y) \Rightarrow x \leq y,$$

and similarly for $x \in \bot_B$. And $g$ is $\mathcal{Y}$-ideal continuous since for $T \in Q$, we get

$$g^{-1}[\mathord{\downarrow}\{T\}] = \{x \in \breve{B} : \mathord{\downarrow}(b,x) \subseteq T\} \cup \bot_B = bT \cup \bot_B \in \mathcal{I}_{\mathcal{Y}}B.$$

Analogously, one defines an $\mathcal{X}$-ideal continuous embedding $h$ of $B$ in $Q$. Furthermore, these two embeddings satisfy for $x \in \breve{A}$ and $y \in \breve{B}$ with $x \bot y$:

$$g(x) \odot h(y) = \overline{t_{\mathcal{X}\mathcal{Y}}}(\mathord{\downarrow}(b,x) \cdot \mathord{\downarrow}(y,b)) = \overline{t_{\mathcal{X}\mathcal{Y}}}(\emptyset) = \emptyset = 0_Q,$$

whereas $(b,b) \in \mathord{\downarrow}(b,x) \cdot \mathord{\downarrow}(y,b) \neq 0_Q$ if $x \bot y$ fails, i.e. $\mathord{\downarrow}x \cap \mathord{\downarrow}y \cap \breve{B} \neq \overline{\emptyset}$. This proves the equivalence $x \bot y \Leftrightarrow g(x) \cdot h(y) = 0_Q$ for $x,y \in \breve{B}$. The cases $x \in \bot_A$ or $y \in \bot_B$ are easy.

(h) $\Rightarrow$ (a): Assume $Y \in \mathcal{Y}$ and $x \bot Y$, i.e. $g(x){\cdot}h(y) = 0_Q$ for all $y \in Y$, hence $g(x) \cdot \bigvee h[Y] = 0_Q$. Since $h$ is $\mathcal{Y}$-ideal continuous, we get $h[\Delta Y] \subseteq \mathord{\downarrow}\bigvee h[Y]$ and so $g(x){\cdot}h(z){=}0_Q$ for all $z \in \Delta Y$, which is equivalent to $x \bot \Delta Y$. Thus, $B$ is $\bot$-$\mathcal{Y}$-distributive, and analogously, $\mathcal{X}$-continuity of $g$ entails that $B$ is $\bot$-$\mathcal{X}$-distributive. $\square$

The case $\mathcal{X} = \mathcal{Y} = \mathcal{P}B$ and

$$t : \mathcal{Q}\breve{B} \to \mathcal{Q}\breve{B},\ R \mapsto \bigcup\{\breve{\Delta}X \times \breve{\Delta}Y : \emptyset \neq X \times Y \subseteq R\},$$

combined with Corollaries 6.1 and 6.2, amounts to

**Theorem 8.2.** *For a complete lattice $B$, the following are equivalent:*

(a) $B$ *is pseudocomplemented.*
(b) $t$ *is a prenucleus on* $\mathcal{Q}\breve{B}$.
(c) $\overline{t}$ *is a nucleus on* $\mathcal{Q}\breve{B}$.
(d) $B \breve{\otimes} B$ *is a quantic quotient of* $\mathcal{Q}\breve{B}$.
(e) $B \breve{\otimes} B \simeq B \otimes B \simeq Gal(B,B)$ *is a quantale.*
(f) $B \breve{\otimes} B \simeq B \otimes B \simeq Gal(B,B)$ *is pseudocomplemented.*
(g) *There exists a (pre)nucleus $j$ on $\mathcal{Q}\breve{B}$ with $t \leq j$ and $j(\emptyset) = \emptyset$.*
(h) *There exist join-preserving order embeddings $g,h$ of $B$ in a quantale $Q$ such that $x \wedge y = 0_B \Leftrightarrow g(x) \cdot h(y) = 0_Q$.*
(i) *For each $\mathcal{X} \subseteq \mathcal{P}B$, $Ant_{\mathcal{X}}(B,B)$ is a quantale.*
(j) *For each $\mathcal{X} \subseteq \mathcal{P}B$, $Ant_{\mathcal{X}}(B,B)$ is pseudocomplemented.*

The mutually inverse isomorphisms
$$T \mapsto \check{T} = T \smallsetminus \overline{\emptyset} \text{ and } T' \mapsto T' \cup \overline{\emptyset}$$
from Corollary 6.1 transport the multiplication $\odot$ on $A_\mathcal{X} \check{\otimes} B_\mathcal{Y}$ to a multiplication $\overline{\odot}$ on $A_\mathcal{X} \otimes B_\mathcal{Y}$ by
$$R \overline{\odot} S = (\check{R} \odot \check{S}) \cup \overline{\emptyset} = \overline{t_{\mathcal{X}\mathcal{Y}}}(\check{R} \cdot \check{S}) \cup \overline{\emptyset}.$$

But, as we have mentioned at the beginning of this section, the resulting multiplication $\overline{\odot}$ heavily differs from the multiplication obtained by taking the full tensor closure of the relation product. Let us add a small but illustrative example.

**Example 8.1.** For a three-element chain $3 = \{0, 1, 2\}$, the truncated tensor product $3 \check{\otimes} 3$ is a non-commutative quantale without unit:

| $\odot$ | $R_0$ | $R_1$ | $R_2$ | $R_3$ | $R_4$ | $R_5$ |
|---|---|---|---|---|---|---|
| $R_0$ | $R_0$ | $R_0$ | $R_0$ | $R_0$ | $R_0$ | $R_0$ |
| $R_1$ | $R_0$ | $R_1$ | $R_2$ | $R_1$ | $R_2$ | $R_2$ |
| $R_2$ | $R_0$ | $R_1$ | $R_2$ | $R_1$ | $R_2$ | $R_2$ |
| $R_3$ | $R_0$ | $R_3$ | $R_5$ | $R_3$ | $R_5$ | $R_5$ |
| $R_4$ | $R_0$ | $R_3$ | $R_5$ | $R_3$ | $R_5$ | $R_5$ |
| $R_5$ | $R_0$ | $R_3$ | $R_5$ | $R_3$ | $R_5$ | $R_5$ |

$R_0 = \emptyset$
$R_1 = \{(1,1)\}$
$R_2 = \{(1,1),(1,2)\}$
$R_3 = \{(1,1),(2,1)\}$
$R_4 = \{(1,1),(1,2),(2,1)\}$
$R_5 = \{(1,1),(1,2),(2,1),(2,2)\}$

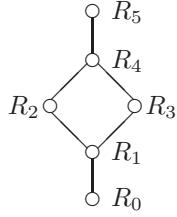

In contrast to that non-trivial quantale, the full tensor closure of the relation product of any two members of $3 \otimes 3$ gives the whole cartesian product $3 \times 3$.

The reader will guess that a result analogous to Theorem 8.1 holds for closure spaces, and that is the case indeed.

**Theorem 8.3.** *A closure space is polarized iff the truncated tensor product $\check{A} \otimes \check{A}$ is a quantale with the multiplication $R \odot S = \overline{R \cdot S}$.*

*Proof.* $A$ is polarized

| | | |
|---|---|---|
| $\Leftrightarrow$ | $\mathcal{C}A$ is pseudocomplemented | (Theorem 7.1) |
| $\Leftrightarrow$ | $\mathcal{C}A \otimes \mathcal{C}A$ is a quantale | (Theorem 8.2) |
| $\Leftrightarrow$ | $A \otimes A$ is a quantale | (Theorem 3.2) |
| $\Leftrightarrow$ | $A \check{\otimes} A$ is a quantale | (Proposition 6.1). $\square$ |

## 9. Unital tensor quantales and boolean algebras

The relation quantales constructed before have a unit (that is, a neutral element for the multiplication) only in a restricted, well decidable situation.

**Proposition 9.1.** *Let $B$ be a poset and $A = At(B)$ the set of its atoms (the minimal elements of $\check{B}$). Then the following conditions are equivalent:*

(a) $B$ *is atomistic, i.e., each element is a join of atoms.*
(b) $I_A = \{(a,a) : a \in A\}$ *is the neutral element of $B \,\check{\otimes}\, B$.*
(c) $B \,\check{\otimes}\, B$ *(respectively, $B \otimes B$) has a neutral element.*

*For bounded $B$, define $i_A \colon B \to B$ by $i_A(0_B) = 1_B$, $i_A|_A = id_A$, $i_A(x) = 0_B$ otherwise. Then the previous conditions are equivalent to the following two:*

(d) *The map $i_A$ is the neutral element of $Gal(B, B)$.*
(e) $Gal(B, B)$ *has a neutral element.*

*Proof.* (a) $\Rightarrow$ (b): $I_A$ is a truncated tensor, being a down-set in $\check{B} \times \check{B}$ and trivially closed under existing joins (notice that $A$ is an antichain). By definition, $I_A \cdot T = \{(a,c) \in T : a \in A\}$ for each $T \in B \,\check{\otimes}\, B$, and consequently the generated truncated tensor $I_A \odot T$ is contained in $T$, too. On the other hand, for each $(b,c) \in T$, we find a nonempty set $X \subseteq A$ with $\bigvee X = b$. Now, $a \in A \cap {\downarrow} b$ entails $(a,c) \in T$ for each $a \in X$, and $(a,a) \in I_A$ yields $(a,c) \in I_A \cdot T$; thus, $X \times \{c\} \subseteq I_A \cdot T \subseteq I_A \odot T$, and therefore $(b,c) \in I_A \odot T$. In all, this proves the equation $I_A \odot T = T$, and $T \odot I_A = T$ is analogous.

(b) $\Rightarrow$ (c): The lattice isomorphism between $B \,\check{\otimes}\, B$ and $B \otimes B$ sends the multiplication $\odot$ to $\overline{\odot}$, hence the unit $I_A$ to the unit $I_A \cup \overline{\emptyset}$.

(c) $\Rightarrow$ (a): Assume $I$ is a neutral element of $B \,\check{\otimes}\, B$. For $(a,b) \in I$, the principal truncated tensor ${\downarrow}(b,a) \subseteq \check{B} \times \check{B}$ is an element of $B \,\check{\otimes}\, B$, whence $(a,a) \in I \odot {\downarrow}(b,a) = {\downarrow}(b,a)$ and $(b,b) \in {\downarrow}(b,a) \odot I = {\downarrow}(b,a)$; thus, $a \leq b$ and $b \leq a$, i.e. $a = b$. Next, assume $0_B < c \leq a$. Then $(a,c) \in I$, and so $a = c$ must be an atom. Now, given any $x \in B$, put $A(x) = \{a \in A : (a,a) \in I, a \leq x\}$, and consider a $y \in B$ with $A(x) \subseteq {\downarrow}y$. For each $(a,b) \in I \cdot {\downarrow}(x,x)$ we have $a \leq y$ (since $(a,a) \in I$) and $b \leq x$; therefore, ${\downarrow}(x,x) = I \odot {\downarrow}(x,x)$ must be contained in the principal truncated tensor ${\downarrow}(y,x)$; it follows that $x \leq y$, and then $x = \bigvee A(x)$.

(b) $\Leftrightarrow$ (d): By Corollary 6.2, $i_A$ is the Galois map $\check{f}_{I_A}$ corresponding to $I_A$.

(c) $\Leftrightarrow$ (e): Use the isomorphisms $B \,\check{\otimes}\, B \simeq B \otimes B \simeq Gal(B, B)$. □

**Lemma 9.1.** *If $B$ and $C$ are atomistic posets then so is $B \,\check{\otimes}\, C \simeq B \otimes C$.*

*Proof.* It is easy to see that $At(B \,\check{\otimes}\, C) = \{{\downarrow}(b,c) : b \in At(B), c \in At(C)\}$. Hence, if $B$ and $C$ are atomistic then for each $T \in B \,\check{\otimes}\, C$, one obtains

$$T = \bigvee \{{\downarrow}(b,c) : b\,T\,c\} = \bigvee \{S \in At(B \,\check{\otimes}\, C) : S \subseteq T\}. \qquad \square$$

Recall that the atomi(sti)c boolean complete lattices (*ABC-lattices*) are nothing but the isomorphic copies of power set lattices ($B$ is *atomic* if each non-zero element is above an atom). From Theorem 8.1 and Proposition 9.1, we now conclude:

**Theorem 9.1.** *For a complete lattice $B$, the following are equivalent:*

(a) $B$ is atomic and boolean (*an ABC-lattice*).
(b) $B$ is atomistic and pseudocomplemented.
(c) $B \,\check{\otimes}\, B \simeq B \otimes B \simeq Gal(B,B)$ is a unital quantale.
(d) $B \,\check{\otimes}\, B$ is isomorphic to the quantale of all relations on a set.

*Proof.* For (b)$\Rightarrow$(a), note that for $a \in At(B)$ and all $b \in B$, one has $a \leq b$ or $a \wedge b = 0_B$, i.e. $a \leq b^*$. Hence, $b \vee b^* = \bigvee At(B) = 1_B$, and by the Glivenko-Frink Theorem [5, 23], this forces $B$ to be Boolean. (For a quick direct proof of distributivity, let $a \in At(B)$, $a \leq b \wedge (c \vee d)$. If $a \not\leq c, d$ then $a \leq c^* \wedge d^* = (c \vee d)^*$, in contrast to $0_B < a \leq c \vee d$. Hence, $a \leq b \wedge c$ or $a \leq b \wedge d$, and as $B$ is a atomistic, this proves $b \wedge (c \vee d) = (b \wedge c) \vee (b \wedge d)$.)
The equivalence (a)$\Leftrightarrow$(d) will follow from Lemma 10.1 and Theorem 10.1. $\square$

In contrast to the *right closure operation* $r = t_{\emptyset \mathcal{P}B}$ and the *left closure operation* $\ell = t_{\mathcal{P}B\, \emptyset}$, the operation $t = t_{\mathcal{P}B\, \mathcal{P}B}$ is frequently only a *preclosure* operation. The lack of idempotency already occurs for finite boolean lattices.

**Example 9.1.** Let $B$ be a boolean algebra with more than four elements. Then there is an $a \in B \smallsetminus \{0,1\}$ whose complement $a^*$ is not an atom, and we find a $b \in B$ with $0 < b < a^*$. Put $c = a^* \wedge b^*$ and $M = \{a,b,c\}$. Then $B$ has the following eight-element boolean subalgebra:

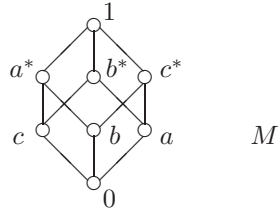

For the lower relation $R = \downarrow\{(x,y) \in M \times M : x \neq y\}$, a quick computation yields $t(R) = \downarrow\{(x,x^*) : x \in M\}$, but $(1,1) \in t^2(R)$.

Our final observation is a bit disappointing:

**Corollary 9.1.** *An unbounded $T_0$ closure space $A$ is discrete if and only if $A \otimes A$ is a unital quantale.*

*Proof.* $A$ is discrete, i.e. $\mathcal{C}A = \mathcal{P}\mathcal{U}A$

$\Leftrightarrow \mathcal{C}A$ is an ABC-lattice $\hspace{2em}$ ($A$ unbounded and $T_0$)
$\Leftrightarrow \mathcal{C}A \otimes \mathcal{C}A$ is a unital quantale $\hspace{2em}$ (Theorem 9.1)
$\Leftrightarrow A \otimes A$ is a unital quantale $\hspace{2em}$ (Theorem 3.2). $\square$

## 10. Semicategories with tensor products as hom-sets

If one wishes to take tensor products of $\perp$-distributive or pseudocomplemented posets or lattices as hom-sets in a category, one has to overcome a little technical problem: while the composition property is saved by the hypothesis of pseudocomplementation (see below), identity morphisms are rather rare, as our results in Section 9 demonstrate. However, Lemma 8.2 assures that we can build at least so-called *semicategories* [24, 27, 35], satisfying the same axioms as categories, with exception of the postulate that each object carries an identity morphism. Such semicategories rather rarely occur in the literature, and if they do, the considerations mostly are confined to *small* semicategories – a restriction we do not want to adopt here. By an obvious adjunction of identity morphisms, one may turn every semicategory into a category, but this may destroy important structural properties; for example, adding a new identity to a (truncated) tensor quantale with the obvious extension of the multiplication gives a monoid, but not a quantale.

The largest semicategory we can construct in the aforementioned way has as objects *partially $\perp$-distributive posets*. These are augmented posets $B_\mathcal{Y} = (B, \mathcal{Y})$ where $B$ is a $\perp$-$\mathcal{Y}$-distributive poset. As hom-set between any two such objects $A_\mathcal{X}$ and $B_\mathcal{Y}$ we take the truncated tensor product $A_\mathcal{X} \check{\otimes} B_\mathcal{Y}$.

**Proposition 10.1.** *The partially $\perp$-distributive posets with the down-set lattices $A_\emptyset \check{\otimes} B_\emptyset$ as hom-sets form a semicategory, and so do the partially $\perp$-distributive posets with the truncated tensor products $A_\mathcal{X} \check{\otimes} B_\mathcal{Y}$ as hom-sets. There is a semifunctor from the former to the latter semicategory, acting identically on objects and sending $R \in A_\emptyset \check{\otimes} B_\emptyset$ to $\overline{t_{\mathcal{X}\mathcal{Y}}}(R) \in A_\mathcal{X} \check{\otimes} B_\mathcal{Y}$.*

*Proof.* The associative law for the composition of lower relations is clear, since it is just the classical relation product; for the morphisms which are members of the truncated tensor products $A_\mathcal{X} \check{\otimes} B_\mathcal{Y}$, the associative law is assured by the homomorphism property of the semifunctor described above, which in turn is an immediate consequence of Lemma 8.2. $\square$

Let us establish an analogous result for closure spaces. Composing the isomorphisms from Theorem 3.2 and Proposition 6.1, we get isomorphisms

$$\mathcal{C}_{AB} : A \otimes B \to \mathcal{C}A \check{\otimes} \mathcal{C}B, \ T \mapsto \{(X, Y) \in \mathcal{C}A \times \mathcal{C}B : X \neq \overline{\emptyset}, Y \neq \overline{\emptyset}, X \times Y \subseteq T\}.$$

**Lemma 10.1.** *For closure spaces $A, B, C$ and $R \in A \otimes B$, $S \in B \otimes C$, $T \in A \otimes C$, one has the equivalence $R \cdot S \subseteq T \Leftrightarrow \mathcal{C}_{AB}(R) \cdot \mathcal{C}_{BC}(S) \subseteq \mathcal{C}_{AC}(T)$ and the equality $\mathcal{C}_{AC}(R \odot S) = \mathcal{C}_{AB}(R) \odot \mathcal{C}_{BC}(S)$.*

*Proof.* Suppose $R \cdot S \subseteq T$. For $(X, Z) \in \mathcal{C}_{AB}(R) \cdot \mathcal{C}_{BC}(S)$ there is a $Y \in \mathcal{C}B$, $Y \neq \overline{\emptyset}$, with $(X, Y) \in \mathcal{C}_{AB}(R)$ and $(Y, Z) \in \mathcal{C}_{BC}(S)$; hence, $X \times Y \subseteq R$, $Y \times Z \subseteq S$, and, as $Y \neq \emptyset$, also $X \times Z \subseteq R \cdot S \subseteq T$ and $(X, Z) \in \mathcal{C}_{AC}(T)$.
Conversely, assume $\mathcal{C}_{AB}(R) \cdot \mathcal{C}_{BC}(S) \subseteq \mathcal{C}_{AC}(T)$ and $(x, z) \in R \cdot S$, say $xRySz$. Then we get $\overline{x} \times \overline{y} \subseteq R$ and $\overline{y} \times \overline{z} \subseteq S$, hence $(\overline{x}, \overline{y}) \in \mathcal{C}_{AB}(R)$, $(\overline{y}, \overline{z}) \in \mathcal{C}_{BC}(S)$, and therefore $(\overline{x}, \overline{z}) \in \mathcal{C}_{AB}(R) \cdot \mathcal{C}_{BC}(S) \subseteq \mathcal{C}_{AC}(T)$, whence $(x, z) \in \overline{x} \times \overline{z} \subseteq R \cdot S$.

Now, by the proven equivalence and the isomorphism property of $\mathcal{C}_{AC}$,

$$\mathcal{C}_{AC}(R \odot S) = \mathcal{C}_{AC}(\bigcap \{T \in A \otimes B : R \cdot S \subseteq T\})$$
$$= \bigcap \{\mathcal{C}_{AC}(T) : \mathcal{C}_{AB}(R) \cdot \mathcal{C}_{BC}(S) \subseteq \mathcal{C}_{AC}(T)\}$$
$$= \bigcap \{\mathcal{T} \in \mathcal{C}A \,\breve{\otimes}\, \mathcal{C}B : \mathcal{C}_{AB}(R) \cdot \mathcal{C}_{BC}(S) \subseteq \mathcal{T}\} = \mathcal{C}_{AB}(R) \odot \mathcal{C}_{BC}(S). \quad \square$$

We are now in a position to define a semifunctor $\mathcal{C}$ from the semicategory **PCS** of polarized closure spaces with tensor products as hom-sets to the semicategory **PCL** of pseudocomplemented complete lattices with truncated tensor products as hom-sets, by taking $\mathcal{C} = \mathcal{C}_{AB}$ on $\mathbf{PCS}(A, B) = A \otimes B$. Then $\mathcal{C}$ preserves the composition by Lemma 10.1, and it induces bijections between the hom-sets $A \otimes B$ and $\mathcal{C}A \,\breve{\otimes}\, \mathcal{C}B$. Furthermore, it is dense, since every pseudocomplemented complete lattice $C$ is isomorphic to the closed set lattice $\mathcal{C}A$ of the polarized cut space $A = (C, \mathcal{N}C)$.

**Theorem 10.1.** *By virtue of the semifunctor $\mathcal{C}$, the semicategories **PCS** and **PCL** are equivalent. Similarly, **PCL** is equivalent to the semicategory of partially $\perp$-distributive posets and truncated tensor products as hom-sets, via the ideal semifunctor sending $A_\mathcal{X}$ to $\mathcal{I}_\mathcal{X} A$ and mapping $A_\mathcal{X} \,\breve{\otimes}\, B_\mathcal{Y}$ to $\mathcal{I}_\mathcal{X} A \,\breve{\otimes}\, \mathcal{I}_\mathcal{Y} B$.*

**Corollary 10.1.** *The pseudocomplemented complete lattices together with the sets of Galois maps as hom-sets form a semicategory whose largest subcategory is that of ABC-lattices, which is equivalent to the category of sets and relations as morphisms.*

In fact, by the distributive laws (which follow from Lemma 2.1), one obtains so-called (*semi*)*quantaloids* as considered by Garraway [24] and Stubbe [35].

Faculty for Mathematics and Physics, Leibniz Univ. Hannover, Germany
*E-mail address*: erne@math.uni-hannover.de

CMUC, Department of Mathematics, Univ. of Coimbra, Portugal
*E-mail address*: picado@mat.uc.pt